\documentclass[reqno,12pt]{amsart}

\pdfoutput=1

\usepackage{color} 
\usepackage{ifpdf}
\ifpdf
    \usepackage[pdftex]{graphicx}
    \usepackage[pdftex]{hyperref}
    \hypersetup{
        unicode=false,          
        pdftoolbar=true,        
        pdfmenubar=true,        
        pdffitwindow=false,     
        pdfstartview={FitH},    
        pdftitle={MCP Article},      
        pdfauthor={Michael Holst},   
        pdfsubject={Mathematics},    
        pdfcreator={Michael Holst},  
        pdfproducer={Michael Holst}, 
        pdfkeywords={PDE, analysis, mathematical physics}, 
        pdfnewwindow=true,      
        colorlinks=true,        
        linkcolor=red,          
        citecolor=blue,         
        filecolor=magenta,      
        urlcolor=cyan           
    }
    
    \def\myfigpdf#1#2{\includegraphics[height=#2]{#1.pdf}}
    \def\myfig#1#2{\includegraphics[height=#2]{#1.png}}
\else
    \usepackage{graphicx}
    \usepackage{pstricks}

    \def\myfig#1#2{\includegraphics[height=#2]{#1.ps}}
\fi

\usepackage{times}
\usepackage{amsfonts}

\usepackage{amsmath}
\usepackage{amsthm}
\usepackage{amssymb}
\usepackage{amsbsy}
\usepackage{amscd}

\usepackage{enumerate}
\usepackage{verbatim}
\usepackage{subfigure}

\newtheorem{theorem}{Theorem}[section]

\newtheorem{definition}[theorem]{Definition}

\newtheorem{algorithm}[theorem]{Algorithm}

\numberwithin{equation}{section}  

  \newcounter{mnote}
  \let\oldmarginpar\marginpar
    \renewcommand\marginpar[1]{\-\oldmarginpar[\raggedleft\footnotesize #1]%
    {\raggedright\footnotesize #1}}

\definecolor{myblue}{rgb}{0.2,0.2,0.7}
\definecolor{mygreen}{rgb}{0,0.6,0}
\definecolor{mycyan}{rgb}{0,0.6,0.6}
\definecolor{myred}{rgb}{0.9,0.2,0.2}
\definecolor{mymagenta}{rgb}{0.9,0.2,0.9}
\definecolor{mywhite}{rgb}{1.0,1.0,1.0}
\definecolor{myblack}{rgb}{0.0,0.0,0.0}

\newcommand{\beq}{\begin{equation}}
\newcommand{\eeq}{\end{equation}}
\newcommand{\beqa}{\begin{eqnarray}}
\newcommand{\eeqa}{\end{eqnarray}}
\def\cal{\mathcal}

\begin{document}

\title[Local Refinement and Multilevel Preconditioning]
      {Local Refinement and Multilevel Preconditioning: \\
       Implementation and Numerical Experiments}

\author[B. Aksoylu]{Burak Aksoylu}
\email{baksoylu@cs.caltech.edu}
\address{Department of Computer Science,
         California Institute of Technology,
         Pasadena, CA 91125, USA}
\thanks{The first author was supported in part by the Burroughs Wellcome
        Fund through the LJIS predoctoral training program at UC San Diego, in
        part by NSF (ACI-9721349, DMS-9872890), and in part by DOE
        (W-7405-ENG-48/B341492). Other support was provided by Intel,
        Microsoft, Alias$|$Wavefront, Pixar, and the Packard Foundation.}

\author[S. Bond]{Stephen Bond}
\email{bond@ucsd.edu}
\address{Department of Chemistry and Biochemistry and
        Department of Mathematics,
        University of California at San Diego,
        La Jolla, CA 92093, USA}
\thanks{The second author was supported in part by the Howard Hughes Medical
        Institute, and in part by NSF and NIH grants to J.~A.~McCammon.}

\author[M. Holst]{Michael Holst}
\email{mholst@math.ucsd.edu}
\address{Department of Mathematics\\
         University of California at San Diego\\ 
         La Jolla, CA 92093, USA}
\thanks{The third author was supported in part by NSF CAREER Award~9875856 and
        in part by a UCSD Hellman Fellowship.}

\date{February 13, 2003}


\begin{abstract}
In this paper, we examine a number of additive and multiplicative
multilevel iterative methods and preconditioners in the setting of
two-dimensional local mesh refinement. While standard multilevel methods are 
effective for uniform refinement-based discretizations of elliptic equations,
they tend to be less effective for algebraic systems
which arise from discretizations on locally refined meshes, 
losing their optimal behavior in both storage and computational 
complexity. Our primary focus here is on BPX-style additive and 
multiplicative multilevel preconditioners, and on various stabilizations 
of the additive and multiplicative hierarchical basis method (HB), and their
use in the local mesh refinement setting. In the first two papers of this 
trilogy, it was shown that both BPX and wavelet stabilizations of HB have 
uniformly bounded conditions numbers on several
classes of locally refined 2D and 3D meshes based on fairly standard
(and easily implementable) red and red-green mesh refinement algorithms.
In this third article of the trilogy, we describe in detail the
implementation of these types
of algorithms, including detailed discussions of the datastructures and
traversal algorithms we employ for obtaining optimal storage and computational
complexity in our implementations.
We show how each of the algorithms can be implemented using standard datatypes
available in languages such as C and FORTRAN, so that the resulting algorithms
have optimal (linear) storage requirements, and so that the resulting
multilevel method or preconditioner can be applied with optimal (linear)
computational costs.
Our implementations are performed in both C and MATLAB using the
Finite Element ToolKit (FEtk), an open source finite element software
package.
We finish the paper with a sequence of numerical experiments illustrating
the effectiveness of a number of BPX and stabilized HB variants for several
examples requiring local refinement.
\end{abstract}

\maketitle


\vspace*{-1.2cm}
{\footnotesize
\tableofcontents
}


\section{Introduction}

While there are a number of effective (often optimal) multilevel methods
for uniform refinement-based discretizations of elliptic equations,
only a handful of these methods are effective for algebraic systems
which arise from discretizations on locally refined meshes, and these
remaining methods are typically suboptimal in both storage
and computational complexity.
In this paper, we examine a number of additive and multiplicative
multilevel iterative methods and preconditioners, specifically for
two-dimensional local mesh refinement scenarios.
Our primary focus is on Bramble, Pasciak, and Xu (BPX)-style additive and 
multiplicative multilevel preconditioners, and on stabilizations of the 
additive and multiplicative hierarchical basis method (HB).
In~\cite{Ak01phd,AkHo02-p1,AkHo02-p2}, it was shown that both BPX and wavelet
stabilizations of HB have uniformly bounded conditions numbers on several
classes of locally refined 2D and 3D meshes based on fairly standard
(and easily implementable) red and red-green mesh refinement algorithms.
In this article, we describe in detail the implementation of these types
of algorithms, including detailed discussions of the datastructures and
traversal algorithms we employ for obtaining optimal storage and computational
complexity in our implementations.
We show how each of the algorithms can be implemented using standard datatypes
available in languages such as C and FORTRAN, so that the resulting algorithms
have optimal (linear) storage requirements, and so that the resulting
multilevel method or preconditioner can be applied with optimal (linear)
computational costs.
Our implementations are performed in both C and MATLAB using the
Finite Element ToolKit (FEtk), an open source finite element software
package.
We also present a sequence of numerical experiments illustrating
the effectiveness of a number of BPX and stabilized HB variants for
examples requiring local refinement.

The problem class of interest for our purposes here 
is linear second order partial
differential equations (PDE) of the form:
\begin{equation}  \label{modelProb}
- \nabla \cdot (p ~ \nabla u) + q ~ u = f,~~~u \in H_0^1(\Omega).
\end{equation}
Here, $f \in L_2(\Omega),~p,~q \in L_\infty(\Omega)$ and
$p:\Omega \rightarrow L(\Re^d, \Re^d),~~ q:\Omega \rightarrow \Re$,
where $p$ is a symmetric positive definite matrix, and $q$ is
nonnegative. Let ${\cal T}_0$ be a shape regular and quasiuniform initial
partition of $\Omega$ into a finite number of $d$-simplices, and
generate ${\cal T}_1, {\cal T}_2, \ldots$ by refining the initial
partition using either red-green or red local refinement strategies in
$d=2$ or $d=3$ spatial dimensions.
Let ${\cal S}_j$ be the simplicial linear $C^0$ finite element (FE)
space corresponding to ${\cal T}_j$ equipped with zero
boundary values. The set of nodal basis functions for ${\cal S}_j$ is
denoted by $\{ \phi_i^{(j)} \}_{i=1}^{N_j}$ where $N_j = \mbox{dim}~ {\cal S}_j$
is equal to the number of interior nodes in ${\cal T}_j$.
Successively refined FE spaces will form the following nested sequence:
\begin{equation} \label{nestedFE}
{\cal S}_0 \subset {\cal S}_1 \subset \ldots \subset {\cal S}_j \subset \ldots
\subset H_0^1(\Omega).
\end{equation}
Although the mesh is nonconforming in the case of red refinement, ${\cal S}_j$
is used within the framework of conforming FE methods for
discretizing~(\ref{modelProb}).

Let the bilinear form and the linear functional representing the weak
formulation of~(\ref{modelProb}) be denoted as
$$
a(u,v) = \int_{\Omega} p~ \nabla u \cdot \nabla v + q~u~v ~dx,
~~~b(v) = \int_{\Omega} f~v~dx,~~~u,v \in H_0^1(\Omega),
$$
and let us consider the following Galerkin formulation: Find
$u \in {\cal S}_j$, such that
\begin{equation} \label{FEproblem}
a(u,v) = b(v),~~~\forall v \in {\cal S}_j.
\end{equation}
Employing the expansion $u=\sum_{i=1}^n u_i^{(j)} \phi_i^{(j)}$ in the
nodal basis for ${\cal S}_j$, problem~(\ref{FEproblem}) reduces to an
algebraic equation of the form:
\begin{equation} \label{ALGproblem}
A^{(j)} {\bf u}^{(j)} = {\bf b}^{(j)} \in \Re^{N_j}
\end{equation}
for the combination coefficients ${\bf u}^{(j)} \in \Re^{N_j}$.
The nodal discretization matrix and vector arise then as:
$$
A^{(j)}_{rs} = a(\phi_s^{(j)},\phi_r^{(j)}),  
~~~~~ {\bf b}^{(j)}_r = b(\phi_r^{(j)}),
 ~~~~ 1 \le r,s \le N_j.
$$

Solving the discretized form of~(\ref{FEproblem}), namely~(\ref{ALGproblem}),
by iterative methods, has been the subject of intensive research because of
the enormous practical impact on a number of application areas in
computational science.
For quality approximation in physical simulation, one is required to use
meshes containing very large numbers of simplices leading to
approximation spaces ${\cal S}_j$ with very large dimension $N_j$.
Only iterative methods which scale well with $N_j$ can be used effectively,
which usually leads to the use of multilevel-type iterative methods
and preconditioners.
Even with the use of such optimal methods for~(\ref{ALGproblem}),
which means methods which scale linearly with $N_j$ in both
memory and computational complexity, the approximation quality requirements
on ${\cal S}_j$ often force $N_j$ to be so large that only
parallel computing techniques can be used to solve~(\ref{ALGproblem}).

To overcome this difficulty one employs adaptive methods, which involves
the use of {\em a posteriori} error estimation to drive {\em local mesh
refinement} algorithms.
This approach leads to approximation spaces ${\cal S}_j$ which are adapted
to the particular target function $u$ of interest, and as a result can
achieve a desired approximation quality with much smaller approximation
space dimension $N_j$ than non-adaptive methods.
One still must solve the algebraic system~(\ref{ALGproblem}), but 
unfortunately most of the available multilevel methods and preconditioners
are no longer optimal, in either memory or computational complexity.
This is due to the fact that in the local refinement setting, the 
approximation spaces $S_j$ do not increase in dimension geometrically
as they do in the uniform refinement setting.
As a result, a single multilevel V-cycle no longer has linear complexity,
and the same difficulty is encountered by other multilevel methods.
Moreover, storage of the discretization matrices and vectors for
each approximation space, required for assembling V-cycle and similar
iterations, no longer has linear memory complexity.

A partial solution to the problem with multilevel methods in the local
refinement setting was provided by the HB method~\cite{BaDu80,BaDuYs88,Ys86}.
This method was based on a direct or hierarchical decomposition of the
approximation spaces ${\cal S}_j$ rather than the overlapping decomposition
employed by the multigrid and BPX method, and therefore by construction had 
linear memory complexity as well as linear computational complexity for a 
single V-cycle-like iteration.
Unfortunately, the HB condition number is not uniformly bounded, leading
to worse than linear overall computational complexity.
While the condition number growth is slow (logarithmic) in two dimensions,
it is quite rapid (geometric) in three dimensions, making it ineffective
in the 3D local refinement setting. Recent alternatives to the HB method, 
including both BPX-like methods~\cite{BrPa93,BPX90} and wavelet-like 
stabilizations of the HB methods~\cite{VaWa3}, provide a final
solution to the condition number growth problem.
It was shown in~\cite{DaKu92} that the BPX preconditioner has
uniformly bounded condition number for certain classes of locally refined
meshes in two dimensions, and more recently in~\cite{AkHo02-p1} it was
shown that the condition number remains uniformly bounded for certain classes
of locally refined meshes in three spatial dimensions.
In~\cite{AkHo02-p2}, it was also shown that wavelet-stabilizations of the
HB method gave rise to uniformly bounded conditions numbers for certain
classes of local mesh refinement in both the two- and three-dimensional
settings.

In view of~\cite{AkHo02-p1} and~\cite{AkHo02-p2}, our interest in this
paper is to examine the practical implementation aspects of both BPX
and stabilized HB iterative methods and preconditioners.
In particular, the remainder of the paper is structured as follows.
In~\S\ref{sec:ml_overview}, we review the algorithms presented
in~\cite{AkHo02-p1} and~\cite{AkHo02-p2}, giving a unified
algorithm framework on which implementations can be based.
The core of the paper is in some sense~\S\ref{sec:implementation},
which describes in detail the datastructures and key algorithms employed
in the implementation of the algorithms.
The focus is on practical realization of optimal (linear) complexity of
the implementations, in both memory and operation complexity.
The FEtk software package which was leveraged for our implementations
is described briefly in~\S\ref{sec:FEtk}.
A sequence of numerical experiments with the implementations
is presented in~\S\ref{sec:numerical}, illustrating the condition number
growth properties of BPX and stabilized HB methods.
Finally, we draw some conclusions in~\S\ref{sec:conc}.

\section{Overview of the multilevel methods}
  \label{sec:ml_overview}

In the first article~\cite{AkHo02-p1} of the trilogy,
it was shown that the BPX preconditioner
was optimal on the meshes under the local 2D and 2D red-green,
as well as local 2D and 3D red, refinement procedures.
The classical BPX preconditioner~\cite{BPX90,XuQi94} 
can be written as an action of the operator $X$ as follows:
\begin{equation} \label{eq:bpxImp}
X u = \sum_{j=0}^J 2^{j(d-2)} \sum_{i=1}^{N_j} (u,\phi_i^{(j)}) \phi_i^{(j)},
~~~ u \in {\cal S}_J.
\end{equation}
Let the prolongation operator from level $j-1$ to $j$ be denoted by $P_{j-1}^j$,
and also denote the prolongation operator from level $j$ to $J$ as:
$$
P_j \equiv P_j^J = P_{J-1}^J  \ldots P_j^{j+1} \in \Re^{N_J \times N_j}, 
$$
where $P_J^J$ is defined to be the identity matrix $I \in \Re^{N_J \times N_J}$.
Then the matrix representation of~(\ref{eq:bpxImp}) becomes~\cite{XuQi94}:
$$
X = \sum_{j=0}^J 2^{j(d-2)} P_j P_j^t.
$$
One can also introduce a version with a smoother $S_j$:
$$
X = \sum_{j=0}^J 2^{j(d-2)} P_j S_j P_j^t.
$$
The preconditioner~(\ref{eq:bpxImp}) can be modified in the hierarchical sense;
\begin{equation} \label{eq:bpxHb}
X_{\mbox{{\tiny {\rm HB}}}} u = 
\sum_{j=0}^J 2^{j(d-2)} \sum_{i=N_{j-1}+1}^{N_j} (u,\phi_i^{(j)}) \phi_i^{(j)},
~~~ u \in {\cal S}_J.
\end{equation}
The new preconditioner corresponds to the additive HB preconditioner in~\cite{Ys86}.
The matrix representation of~(\ref{eq:bpxHb}) is formed from matrices $H_j$
which are simply the tails of the $P_j$ corresponding to newly introduced
degrees of freedom (DOF) in the fine space.
In other words, $H_j \in \Re^{N_J \times (N_j - N_{j-1})}$
is given by only keeping the fine columns (the last $N_j - N_{j-1}$ columns of $P_j$). Hence,
the matrix representation of~(\ref{eq:bpxHb}) becomes:
$$
X_{\mbox{{\tiny {\rm HB}}}} = \sum_{j=0}^J 2^{j(d-2)} H_j H_j^t.
$$

Only in the presence of a geometric increase in the number of DOF, 
the same assumption for optimality of a single classical 
multigrid or BPX iteration, does the cost per iteration remain optimal.
In the case of local refinement, the BPX preconditioner (\ref{eq:bpxImp})
(usually known as additive multigrid)  can easily be suboptimal because 
of the suboptimal cost per iteration (see 
Figure~\ref{figure:flopCountMulAddExpII}). On the other hand, the HB 
preconditioner (\ref{eq:bpxHb}) suffers from a suboptimal iteration count.
The above deficiencies of the preconditioners (\ref{eq:bpxImp}) and
(\ref{eq:bpxHb}) can be overcome by restricting the sum over $i$ in 
(\ref{eq:bpxImp}) only to those nodal basis functions with supports that
intersect the refinement region \cite{BoYs93,BrPa93,DaKu92,Os94book}. 
We call this set {\em 1-ring} of fine DOF,
namely, the set which contains fine DOF and their immediate neighboring coarse
DOF. The following is referred as the BPX preconditioner for local refinement.
\begin{equation} \label{eq:bpxLocal}
X u = 
\sum_{j=0}^J 2^{j(d-2)} 
\sum_{ i \in \mbox{{\tiny {\rm ONERING}}} }
(u,\phi_i^{(j)}) \phi_i^{(j)},
~~~ u \in {\cal S}_J,
\end{equation}
where {\rm ONERING}$=\{1-ring(ii):~ii=N_{j-1}+1, \ldots, {N_j} \}$.

\begin{figure}[ht]
\begin{center}
\mbox{\myfig{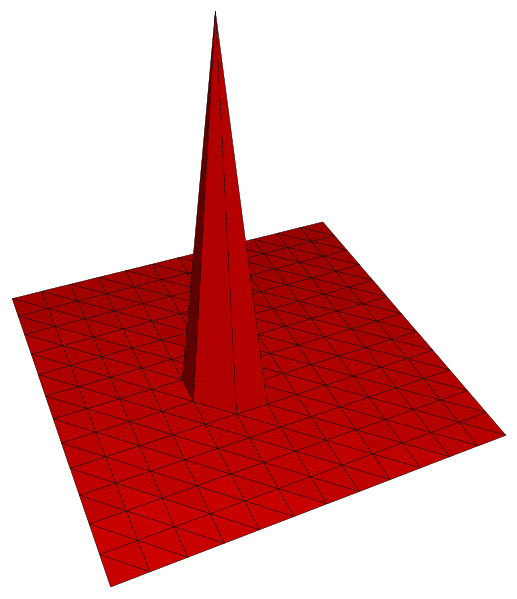}{5cm}} 
\end{center}
\caption{Hierarchical basis function without modification.}
\label{figure:nodalBasis}
\end{figure}

\begin{figure}[ht]
\begin{center}
\mbox{\myfig{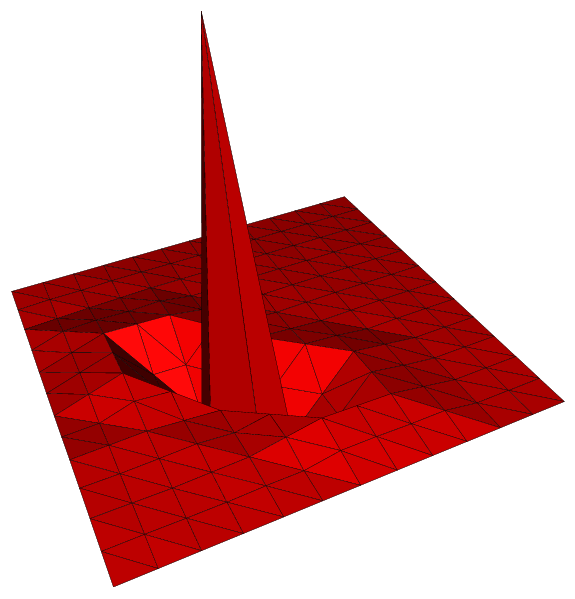}{5cm}} 
\hspace*{1.0cm}
\mbox{\myfig{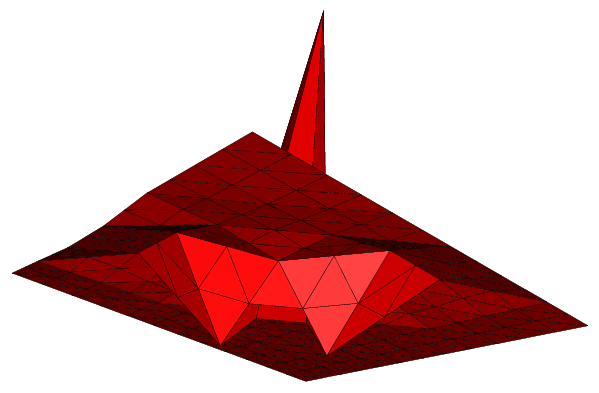}{4.5cm}} 
\end{center}
\caption{Wavelet modified hierarchical basis function with one iteration of
symmetric Gauss-Seidel approximation, upper and lower view.}
\label{figure:basisFuncGS}
\end{figure}

\begin{figure}[ht]
\begin{center}
\mbox{\myfig{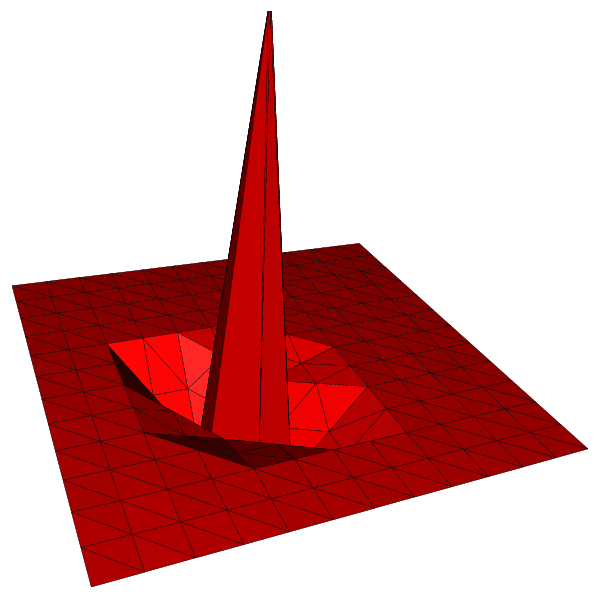}{5cm}} 
\hspace*{1.0cm}
\mbox{\myfig{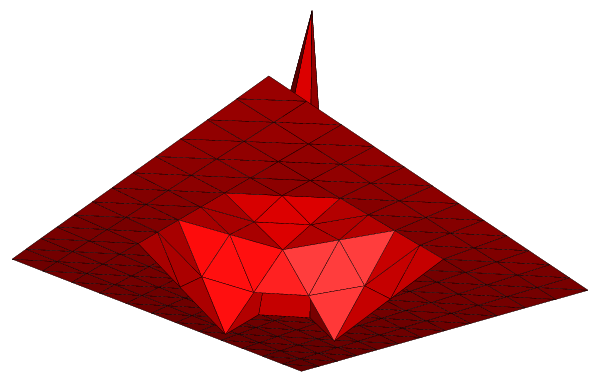}{4.5cm}} 
\end{center}
\caption{Wavelet modified hierarchical basis function with one iteration of
Jacobi approximation, upper and lower view.}
\label{figure:basisFuncJac}
\end{figure}

The BPX decomposition gives rise to basis functions which are not locally 
supported, but they decay rapidly outside a local support region.
This allows for locally supported {\em approximations} as illustrated
in Figures~\ref{figure:nodalBasis}, \ref{figure:basisFuncGS},
and~\ref{figure:basisFuncJac}.

The {\em wavelet modified hierarchical basis
(WMHB) methods}~\cite{VaWa2,VaWa1,VaWa3} can be viewed as an approximation
of the wavelet basis stemming from the BPX decomposition~\cite{Ja92}.
A similar wavelet-like multilevel decomposition approach was taken 
in~\cite{St97}, where the orthogonal decomposition is formed by a discrete 
$L_2$-equivalent inner product. This approach utilizes the same BPX two-level 
decomposition~\cite{St95,St97}.

For adaptive regimes, other primary method of interest is the WMHB method.
The WMHB methods can be described as additive or multiplicative Schwarz
methods. 
In one of the previous papers~\cite{AkHo02-p2} of this trilogy,
it was shown that the additive version of the WMHB method was optimal
under certain types of red-green mesh refinement.
Following the notational framework in~\cite{AkHo02-p2,VaWa3},
this method is defined recursively as follows:
\begin{definition} \label{defn:Dj}
The additive WMHB method $D^{(j)}$ is defined for $j=1,\ldots,J$ as
$$ D^{(j)} \equiv \left[
\begin{array}{ll}
D^{(j-1)}    & 0 \\
0            & B^{(j)}_{22}
\end{array} \right],
$$
with $D^{(0)}=A^{(0)}$.
\end{definition}

With smooth PDE coefficients, optimal results were also established for the
multiplicative version of the WMHB method in~\cite{AkHo02-p2}. Our numerical
experiments demonstrate such optimal results. This method can be written
recursively as:
\begin{definition} \label{defn:Bj}
The  multiplicative WMHB method $B^{(j)}$ is defined as
$$ B^{(j)} \equiv \left[
\begin{array}{ll}
B^{(j-1)}    & A^{(j)}_{12} \\
0            & B^{(j)}_{22}
\end{array} \right]
\left[ \begin{array}{ll}
I                            & 0 \\
B^{(j)^{-1}}_{22} A_{21}^{(j)} & I
\end{array} \right] =
\left[ \begin{array}{ll}
B^{(j-1)}+ A^{(j)}_{12}B^{(j)^{-1}}_{22} A_{21}^{(j)} & A^{(j)}_{12} \\
A_{21}^{(j)} & B^{(j)}_{22}
\end{array} \right],
$$
with $B^{(0)}=A^{(0)}$.
\end{definition}

$A^{(j)}_{12}, A^{(j)}_{21}, A^{(j)}_{22}$ represent subblocks of $A^{(j)}$ and
they correspond to coarse-fine, fine-coarse, and fine-fine interactions of 
DOF at level $j$, respectively. $B^{(j)}_{22}$ denotes an approximation of 
$A^{(j)}_{22}$, e.g. Gauss-Seidel or Jacobi approximation.
For a more complete description of these and related algorithms, 
see~\cite{AkHo02-p1,AkHo02-p2}.

\section{Implementation}
  \label{sec:implementation}

The overall utility of any finite element code depends strongly
on efficient implementation of its core algorithms and data structures.
Theoretical results involving complexity are of little practical importance
if the methods cannot be implemented.  For algorithms involving data 
structures, this usually means striking a balance between storage costs and 
computational complexity.  Finding a minimal representation for a data set is 
only useful if the information can be accessed efficiently. 
\subsection{Sparse Matrix Structures}
Our implementation relies on a total of four distinct sparse 
matrix data structures:  compressed column (COL), compressed row (ROW),
diagonal-row-column (DRC), and orthogonal-linked list (XLN).
Each of these storage schemes attempts to record the location and value of
the nonzeros using a minimal amount of information.  The schemes differ in
the exact representation which effects the speed and manner with which the 
data can be retrieved.  To illustrate how each of these data structures works 
in practice, we consider storing the following sparse matrix:
\begin{equation}
\left[ \begin{array}{ccccc} 1 & 2 &  &  & \\
3 & 4 & 5 &  & 6 \\  & 7 & 8 &   & \\
  &  &  & 9 & 10 \\  & 11 & & 12 & 13
\end{array}
\right] 
\end{equation}

$\bullet$ COL:  The compressed column format is the most commonly 
used sparse matrix type in the literature.  It is the format chosen for
the Harwell-Boeing matrix collection~\cite{DuGrLe89}, and is used in 
production codes such as SuperLU~\cite{DeEiGiLiLiu99}.  In this 
data structure, the nonzeros are 
arranged by column
in a single double-precision array:
\[ \mathrm{A}_{\mathrm{COL}} = \left[ 1, \, 3, \, 2, \, 4, \, 7, \, 11, \, 5, \, 8, \, 9, \, 12, \, 6, \, 10, \, 13 \right]. \]
The indices of $\mathrm{A}$ (often referred to as pointers) corresponding to 
the first entry in each column is then stored in an integer array:
\[ \mathrm{IA}_{\mathrm{COL}} = \left[ 1, \, 3, \, 7, \, 9, \, 11, \, 14 \right].\]
The length of the array $\mathrm{IA}$ is always one greater than the number 
of columns, with the last entry is equal to the number of nonzeros plus one.  
The difference in successive entries in the $\mathrm{IA}$ array reflects the 
number of nonzeros in each column.  If a column has no nonzeros, the index 
from the next column is repeated.  To determine the location of each nonzero 
within its column, the row index of each entry is stored in an integer array:
\[ \mathrm{JA}_{\mathrm{COL}} = \left[ 1, \, 2, \, 1, \, 2, \, 3, \, 5, \, 2, \, 3, \, 4, \, 5, \, 2, \, 4, \, 5 \right]. \]
There is no restriction that the entries are ordered within each column,
only that the columns are ordered.  The memory required to store this
datastructure is: $(nZ + nC + 1) * \mathrm{size}\left(\mathrm{int}\right) + nZ * \mathrm{size}\left(\mathrm{double}\right)$, where $nZ$ and $nC$ are the
number of nonzeros and columns respectively.

$\bullet$ ROW:  The compressed row data structure is just the transpose of the 
compressed column data structure, where the nonzero entries, row pointers,
and column indices are stored in $\mathrm{A}$, $\mathrm{IA}$, and 
$\mathrm{JA}$ respectively:
\[ \mathrm{A}_{\mathrm{ROW}} = \left[ 1, \, 2, \, 3, \, 4, \, 5, \, 6, \, 7, \, 8, \, 9, \, 10, \, 11, \, 12, \, 13 \right], \]
\[ \mathrm{IA}_{\mathrm{ROW}} = \left[ 1, \, 3, \, 7, \, 9, \, 11, \, 14 \right], \quad 
\mathrm{JA}_{\mathrm{ROW}} = \left[ 1, \, 2, \, 1, \, 2, \, 3, \, 5, \, 2, \, 3, \, 4, \, 5, \, 2, \, 4, \, 5 \right]. \]
One should note that since in our example the matrix is {\em structurally symmetric}, the $\mathrm{IA}$ and $\mathrm{JA}$ arrays are identical in both the
ROW and COL cases. The memory required to store this datastructure is: 
$(nZ + nR + 1) * \mathrm{size}\left(\mathrm{int}\right) + nZ * \mathrm{size}\left(\mathrm{double}\right)$, where $nR$ is the number of rows.

$\bullet$ DRC:  The diagonal-row-column format is a structurally symmetric data
structure, which is only valid for square matrices.  In this format, the
diagonal is stored in its own full vector, while the strictly upper and lower
triangular portions are stored in ROW and COL formats respectively.
Leveraging the symmetry in the nonzero structure, the same $\mathrm{IA}$ and 
$\mathrm{JA}$ arrays can be used for the upper and lower triangular parts:
\[ \mathrm{AD}_{\mathrm{DRC}} = \left[ 1, \, 4, \, 8, \, 9, \, 13 \right], 
\] 
\[ \mathrm{AU}_{\mathrm{DRC}} = \left[ 2, \, 5, \, 6, \, 10 \right], \quad
\mathrm{AL}_{\mathrm{DRC}} = \left[ 3, \, 7, \, 11, \, 12 \right]
\]
\[ \mathrm{IA}_{\mathrm{DRC}} = \left[ 1, \, 2, \, 4, \, 4, \, 5, \, 5 \right], \quad 
\mathrm{JA}_{\mathrm{DRC}} = \left[ 2, \, 3, \, 5, \, 5 \right]. \]
The memory required to store this datastructure is less than ROW or COL if
the diagonal is full, and the matrix is structurally symmetric.

\begin{figure}[ht]
\begin{center}
\mbox{\myfigpdf{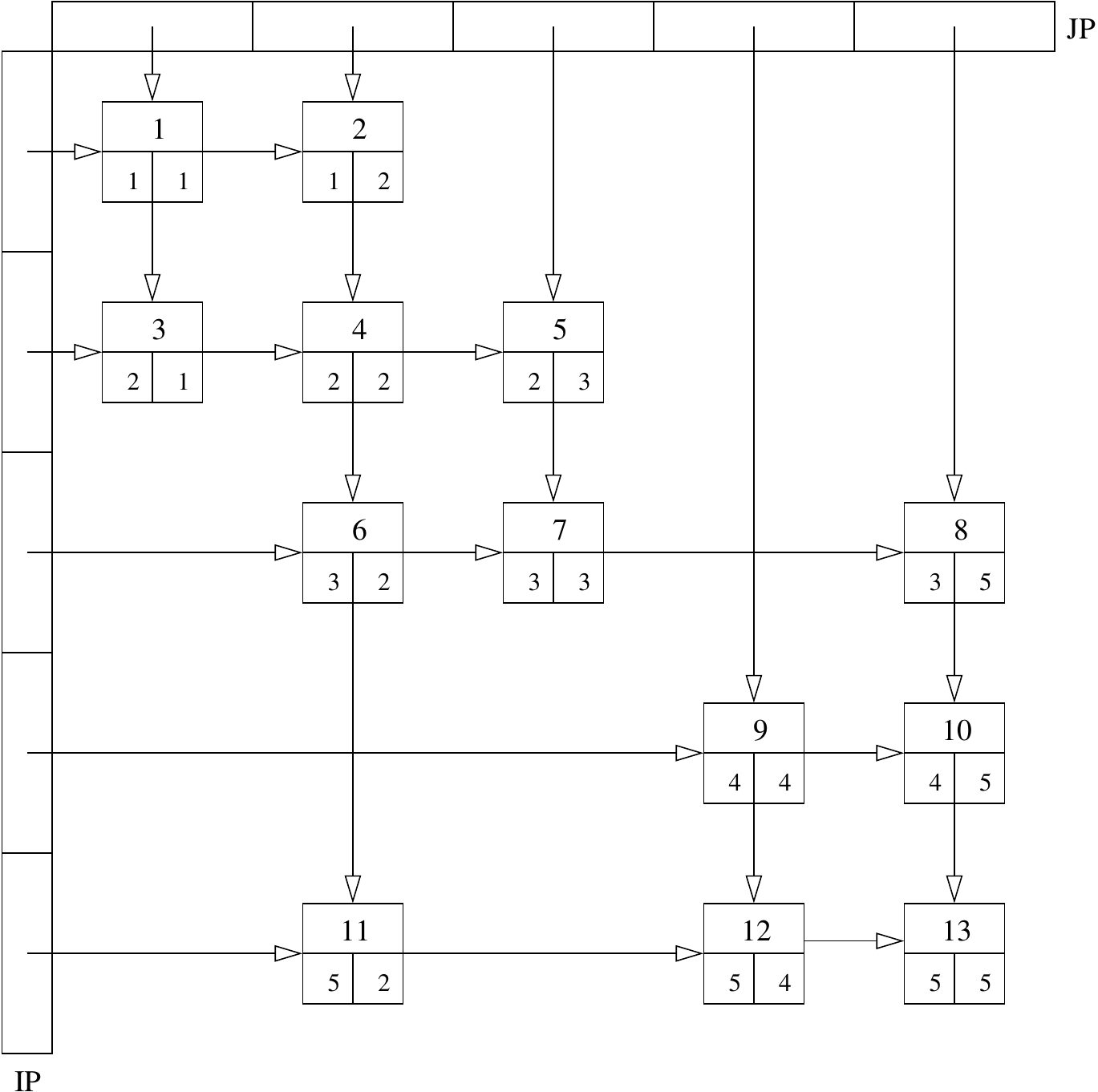}{4.0in}}
\end{center}
\caption{\label{fig:XLN}
An illustration of the XLN datastructure.} 
\end{figure}

$\bullet$ XLN:  The orthogonal-linked list format is the only dynamically 
``fillable'' datastructure used by our methods.  By using variable length 
linked lists, rather than a fixed length array, it is suitable for situations 
where the total number of nonzeros is not known {\it a priori}.
The XLN datastructure is illustrated graphically in Figure~\ref{fig:XLN}.
For each nonzero, there is a link containing the value, row index, column 
index, and pointers to the next in the row and column.  To keep track of
the first link in each row and column, there are two additional pointer 
arrays, $\mathrm{IP}$ and $\mathrm{JP}$.  As long as there are ``order-one''
nonzeros per row, accessing any entry can be accomplished in ``order-one''
time.  The structure can be traversed both rowwise, and columnwise depending
on the situation.  If the matrix is symmetric, only the lower triangular
portion is stored.  The total storage overhead for this structure is:
$nZ * \left( \mathrm{size}\left(\mathrm{double}\right) + 2 * \mathrm{size}\left(\mathrm{int}\right) \right) + \left( nC + nR + 2 \, nZ \right) * \mathrm{size}\left(\mathrm{ptr}\right)$.  Although this is considerably more than the
other three datastructures, one should note that the asymptotic complexity
is still linear in the number of nonzeros.

\subsection{Sparse Matrix Products}
The key preprocessing step in the hierarchical basis methods, is converting
the ``nodal'' matrices and vectors into the hierarchical basis.  This 
operation involves sparse matrix-vector and matrix-matrix products for
each level of refinement.  To ensure that this entire operation has linear
cost, with respect to the number of unknowns, the per-level change of 
basis operations must have a cost of ${\cal O}\left( n_j \right)$, where
$n_j := N_j - N_{j-1} $ is the number of ``new'' nodes on level $j$.  For 
the traditional multigrid algorithm this is not possible, since enforcing 
the variational conditions operates on {\em all} the nodes on each level, 
not just the newly introduced nodes.

The linear operator which converts from the nodal to the hierarchical basis
can be written in terms of a change of basis matrix:
\[ G = \left[ \begin{array}{cc} I & K_{12} \\ K_{21} & I + K_{22} 
\end{array} \right], \]
where $G \in \Re^{N_j \times N_j}$, $K_{12} \in \Re^{N_{j-1} \times n_j}$, 
$K_{21} \in \Re^{n_j \times N_{j-1}}$, and $K_{22} \in \Re^{n_j \times n_j}$.
In this representation, we have assumed that the nodes are ordered with the
nodes $N_{j-1}$ inherited from the previous level listed first, and the $n_j$
new DOF listed second.  For both wavelet modified (WMHB) and 
unmodified hierarchical basis (HB), the $K_{21}$ block represents the 
last $n_j$ rows of the prolongation matrix, $P_{j-1}^{j}$.  In the HB
case, the $K_{12}$ and $K_{22}$ blocks are zero resulting in a very
simple form:
\begin{equation}
G_{\mathrm{hb}} = \left[ \begin{array}{cc} I & 0 \\ K_{21} & I  \end{array} \right]
\end{equation}
For WMHB, the $K_{12}$ and $K_{22}$ blocks are computed using the mass
matrix, which results in the following formula:
\begin{equation}
G_{\mathrm{wmhb}} = \left[ \begin{array}{cc} I &  - \mathrm{inv}\left[ M_{11}^{\mathrm{hb}} \right] M_{12}^{\mathrm{hb}}  \\ 
K_{21} & I   - K_{21} \mathrm{inv}\left[ M_{11}^{\mathrm{hb}} \right] M_{12}^{\mathrm{hb}} \end{array} \right],
\end{equation}
where the $\mathrm{inv}\left[ \cdot \right]$ is some approximation to the
inverse which preserves the complexity.
For example, it could be as simple as the inverse of the diagonal, or a 
low-order matrix polynomial approximation.  The $M^{\mathrm{hb}}$ blocks are 
taken from the mass matrix in the HB basis:
\begin{equation}
M^{\mathrm{hb}} = G_{\mathrm{hb}}^T M^{\mathrm{nodal}} G_{\mathrm{hb}}.
\end{equation}
For the remainder of this section, we restrict our attention to the WMHB
case.  The HB case follows trivially with the two additional subblocks
of $K$ set to zero.  

To reformulate the nodal matrix representation of the bilinear form in 
terms of the hierarchical basis, we must perform a triple matrix product
of the form:
\begin{eqnarray*}  
A_{(j)}^{\mathrm{wmhb}} & = & G_{(j)}^T A_{(j)}^{\mathrm{nodal}} G_{(j)} \\
 & =  & \left( I + K_{(j)}^T \right) A_{(j)}^{\mathrm{nodal}} 
\left( I + K_{(j)} \right).
\end{eqnarray*}
In order to keep linear complexity, we can only copy $A^{\mathrm{nodal}}$
a fixed number of times, i.e. it cannot be copied on every level.  Fixed
size data-structures are unsuitable for storing the product, since predicting 
the nonzero structure of $A_{(j)}^{\mathrm{wmhb}}$ is just as difficult
as actually computing it.  It is for these reasons that we have chosen
the following strategy:  First, copy  $A^{\mathrm{nodal}}$ on the 
finest level, storing the result in an XLN which will eventually become $A^{\mathrm{wmhb}}$.  Second, form the product pairwise, contributing the result to 
the XLN.  Third, the last $n_j$ columns and rows of $A^{\mathrm{wmhb}}$ are 
stripped off, stored in fixed size blocks, and the operation is repeated on 
the next level, using the $A_{11}$ block as the new  $A^{\mathrm{nodal}}$:
\begin{algorithm}
   \label{alg:WMHBSTAB}
(Wavelet Modified Hierarchical Change of Basis)
{\small \begin{itemize}
\item Copy $A_{J}^{\mathrm{nodal}} \rightarrow A^{\mathrm{wmhb}}$ in XLN 
format.
\item While $j > 0$
    \begin{enumerate}
    \item Multiply $ A^{\mathrm{wmhb}} = A^{\mathrm{wmhb}} G $ as 
	\[ \left[ \begin{array}{cc} A_{11} & A_{12} \\
	A_{21} & A_{22} \end{array} \right] += 
\left[ \begin{array}{cc} A_{11} & A_{12} \\
	A_{21} & A_{22} \end{array} \right]
\left[ \begin{array}{cc} 0 & K_{12} \\
	K_{21} & K_{22} \end{array} \right]
\] 

    \item Multiply $ A^{\mathrm{wmhb}} = G^T A^{\mathrm{wmhb}} $ as 
	\[ \left[ \begin{array}{cc} A_{11} & A_{12} \\
	A_{21} & A_{22} \end{array} \right] += 
\left[ \begin{array}{cc} 0 & K_{21}^T \\
	K_{12}^T & K_{22}^T \end{array} \right]
\left[ \begin{array}{cc} A_{11} & A_{12} \\
	A_{21} & A_{22} \end{array} \right]
\] 
    \item Remove $A^{(j)}_{21}$, $A^{(j)}_{12}$, $A^{(j)}_{22}$ blocks of
$A^{\mathrm{wmhb}}$ storing in ROW, COL, and DRC formats respectively.
    \item After the removal, all that remains of $A^{\mathrm{wmhb}}$ is its $A^{(j)}_{11}$ block.
    \item Let j = j - 1, descending to level $j-1$.
    \end{enumerate}
 \item End While.
    \item Store the last $A^{\mathrm{wmhb}}$ as $A_{\mathrm{coarse}}$
\end{itemize} }
\end{algorithm}

We should note that in order to preserve the complexity of the overall
algorithm, all of the matrix-matrix algorithms must be carefully 
implemented.  For example, the change of basis involves computing the
products of $A_{11}$ with $K_{12}$ and $K_{12}^T$.  To preserve storage 
complexity, $K_{12}$ must be kept in compressed column format, COL.
For the actual product, the loop over the columns of $K_{12}$ must be
ordered first, then a loop over the nonzeros in each column, then a loop
over the corresponding row or column in $A_{11}$.
It is exactly for this reason, that one must
be able to traverse $A_{11}$ both by row and by column, which is why we
have chosen an orthogonal-linked matrix structure for $A$ during the 
change of basis (and hence $A_{11}$).  

To derive optimal complexity algorithms for the other products, it is
enough to ensure that the outer loop is always over a dimension
of size $n_j$.  Due to the limited ways in which a sparse matrix 
can be traversed, the ordering of the remaining loops will usually
be completely determined.
Further gains can be obtained in the symmetric case, since
only the upper or lower portion of the matrix needs to be explicitly
computed and stored.

\subsection{The Finite Element ToolKit (FEtk)}
  \label{sec:FEtk}

A number of variations of the methods described above have been
implemented using the Finite Element ToolKit (FEtk)~\cite{Hols2001a}.
FEtk is an open source finite element modeling package which
has been developed by the Holst research group over several years
at Caltech and UC San Diego, with generous contributions from a number
of colleagues.
FEtk consists of a low-level portability library called
MALOC (Minimal Abstraction Layer for Object-oriented C), on top of
which is build a general finite element modeling kernel called
MC (Manifold Code).
Most of the images appearing later in this paper were produced using 
another component of FEtk call SG (Socket Graphics), which is also
built on top of MALOC.
FEtk also includes a fully functional MATLAB version of MC called MCLite,
which shares with MC its datastructures, {\em a posteriori} error estimation
and mesh refinement algorithms, and iterative solution methods.
All of the preconditioners employed in this paper have been implemented by
the authors as ANSI-C class library extensions to MC, and as MATLAB
toolkit-like extensions to MCLite.
The two implementations are mathematically equivalent, although the MCLite 
implementation is restricted to two spatial dimensions.
(The MC-based implementation is both two- and three-dimensional.)
The extensions to MC are distributed as the {\em MCX} library,
and as MATLAB extensions to MCLite are distributed as {\em MCLiteX}.

MALOC, SG, MC, and MCLite are freely redistributable under the
GNU General Public License (GPL).
More information about FEtk can be found at:
\begin{center}
{\tt http://www.fetk.org} 
\end{center}

\section{Numerical Experiments}
  \label{sec:numerical}

\begin{figure}[ht]
\begin{center}
\mbox{
\myfigpdf{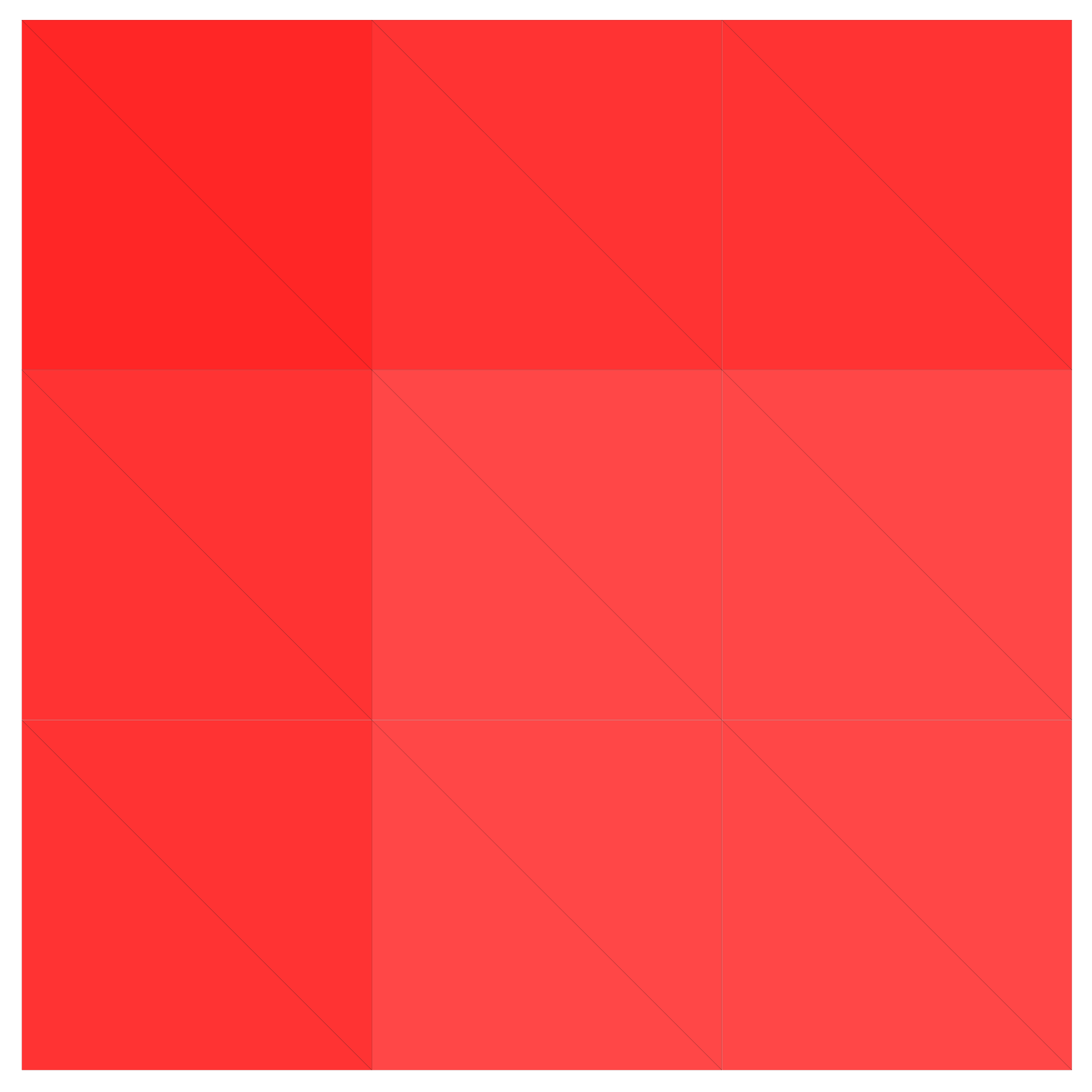}{4.0cm} \hspace{1cm} 
\myfigpdf{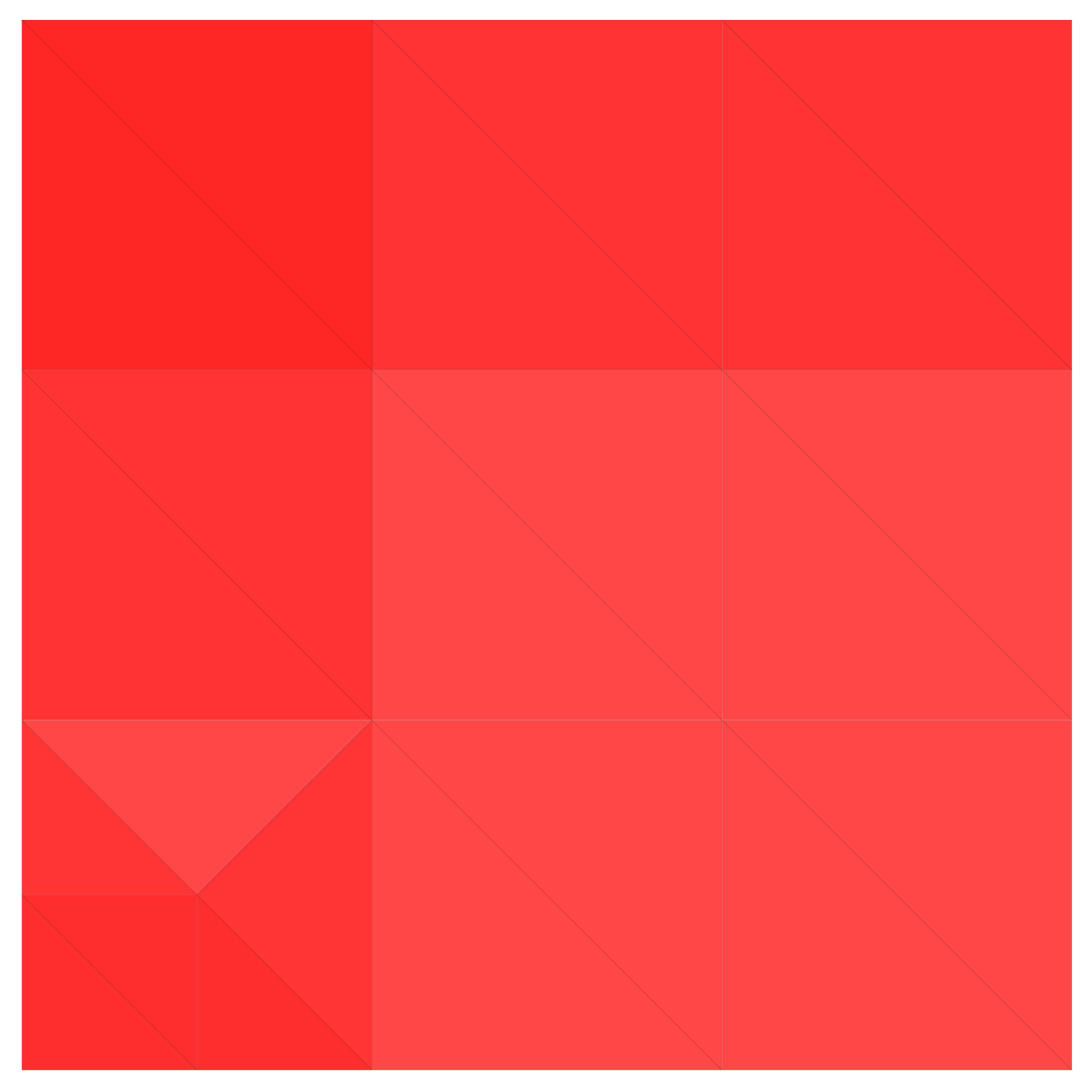}{4.0cm} \hspace{1cm} 
\myfigpdf{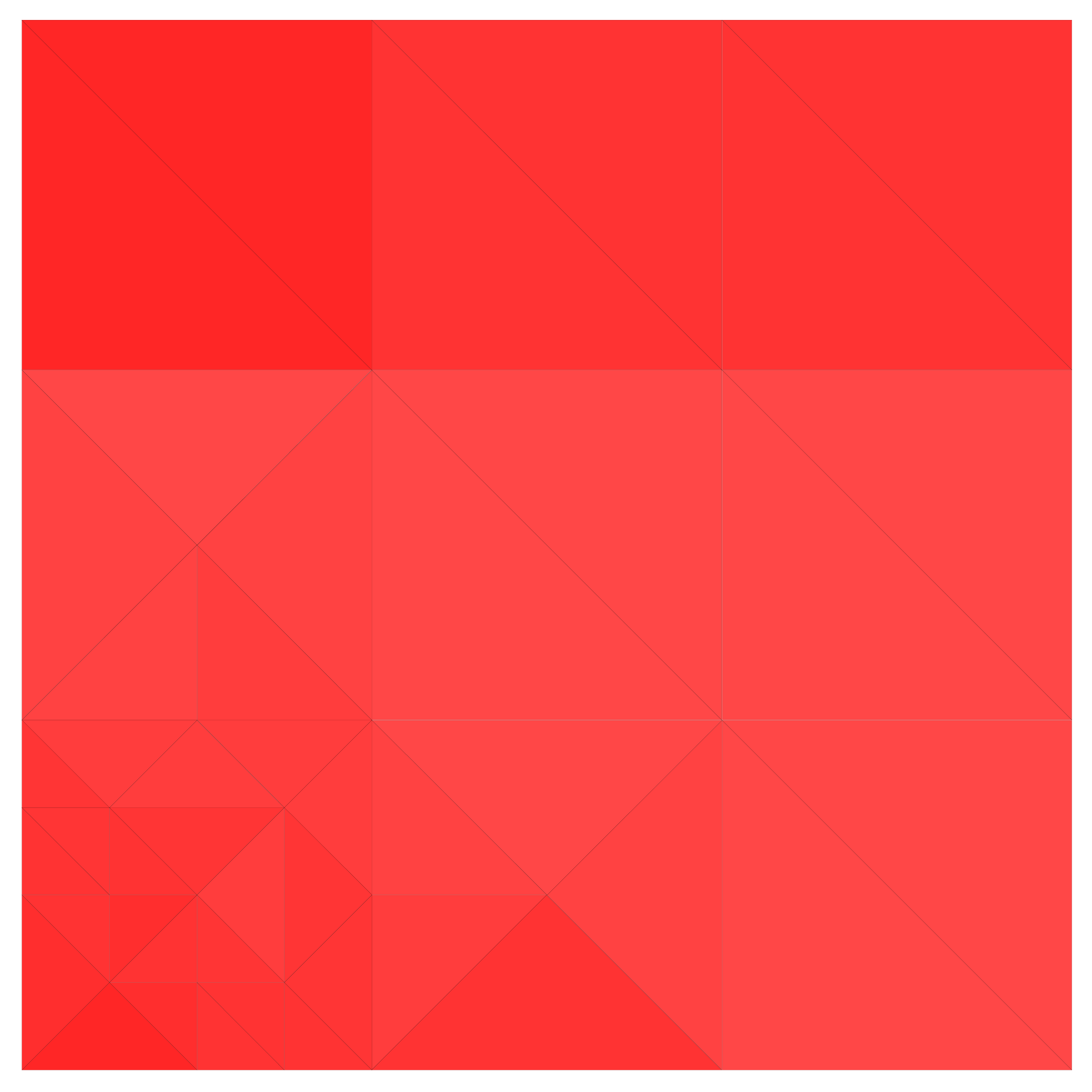}{4.0cm}  
} \\
\mbox{ 
\myfigpdf{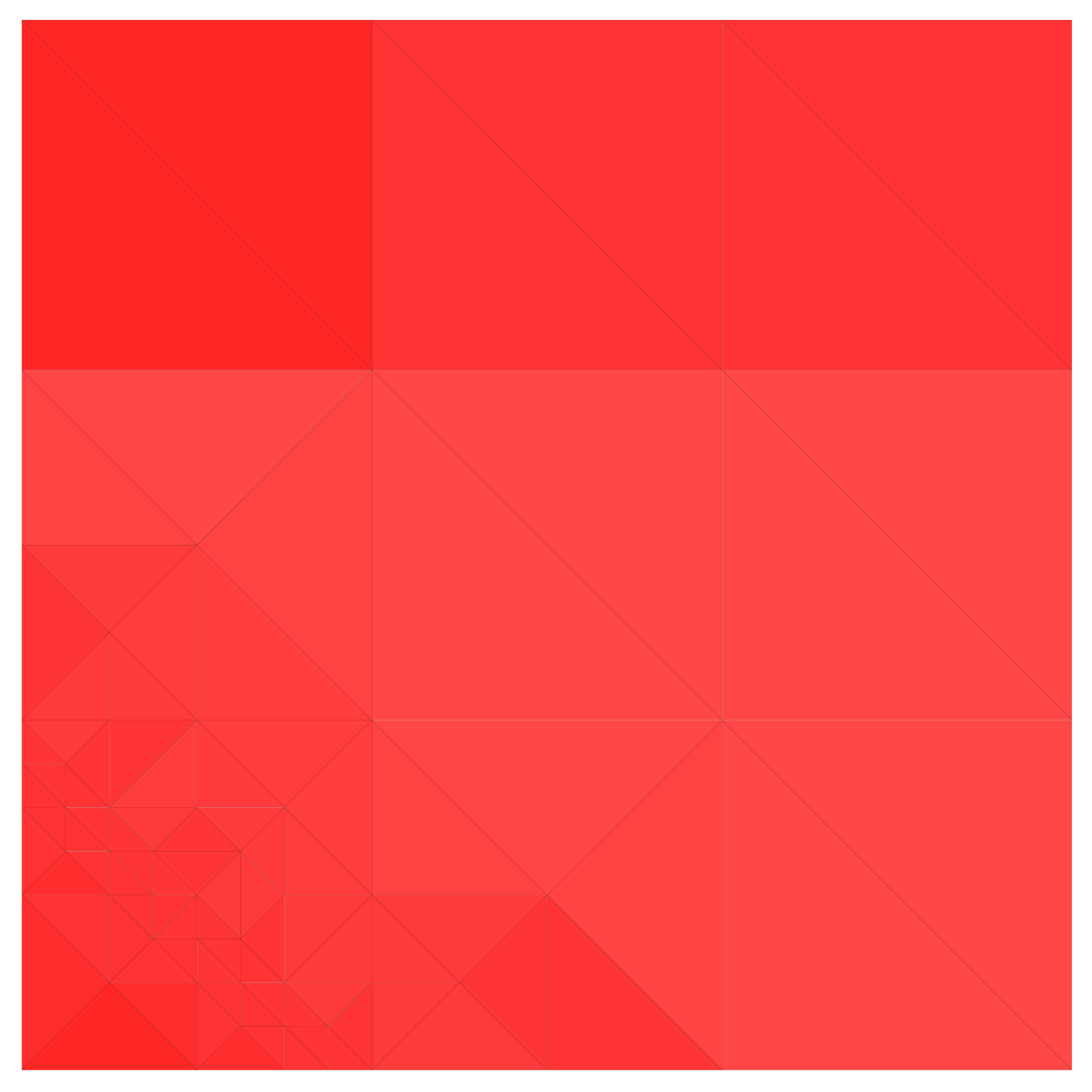}{4.0cm} \hspace{1cm} 
\myfigpdf{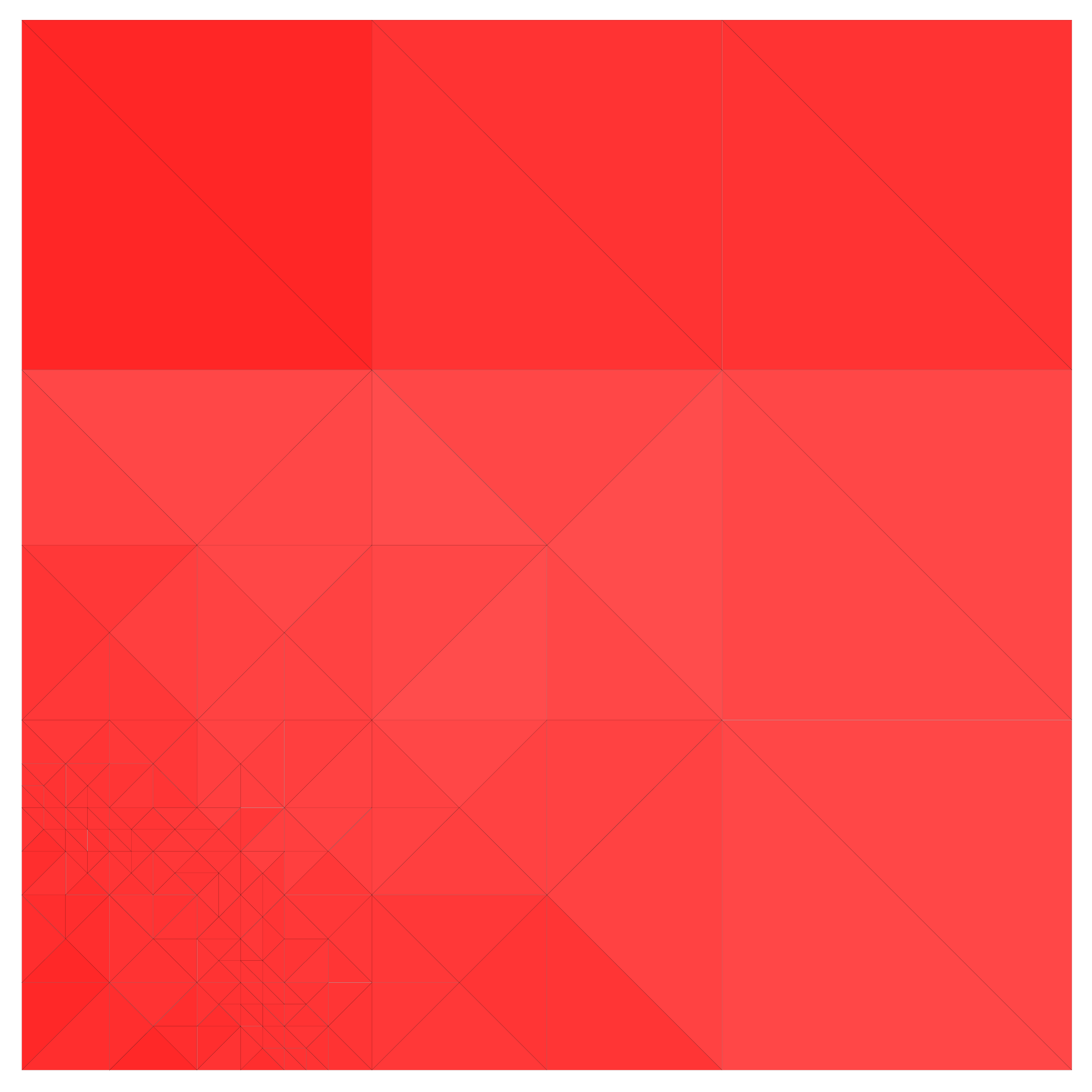}{4.0cm} \hspace{1cm} 
\myfigpdf{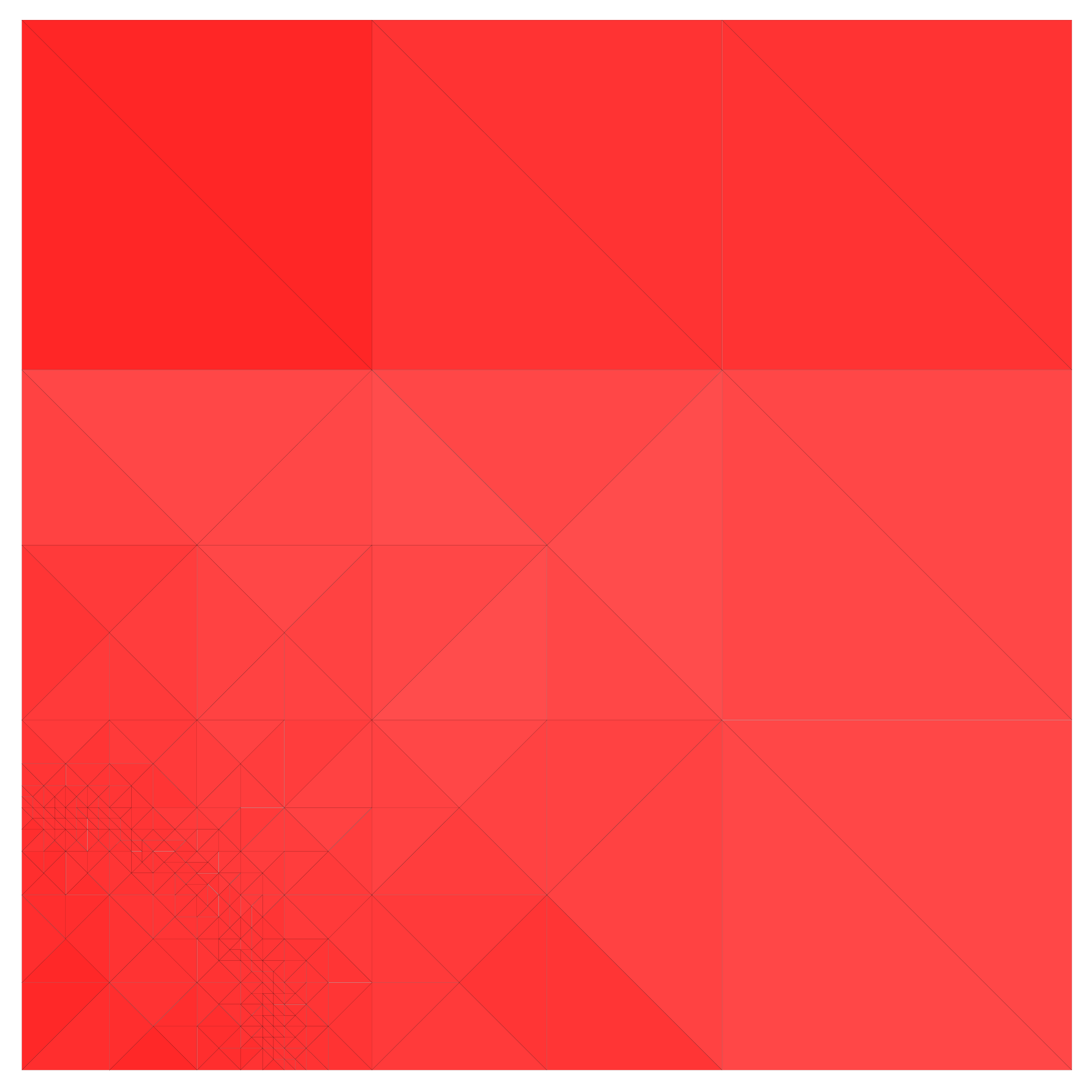}{4.0cm}  
}
\caption{Adaptive mesh, experiment set I.}
\label{figure:meshMclabSpecial}
\end{center}
\end{figure}

Test problem is as follows:
\begin{eqnarray*}
- \nabla \cdot (p~ \nabla u) + q~ u &=& f, 
   ~~~ x \in \Omega \subset \Re^2, \\
n \cdot (p~ \nabla u) &=& g, ~~~ \mbox{on}~ \Gamma_N, \\
  u &=& 0, ~~~ \mbox{on}~ \Gamma_D,
\end{eqnarray*}
where $\Omega = [0,1] \times [0,1]$ and
\[
	p = \left [ \begin{array}{cc} 1 & 0 \\ 0 & 1 \end{array} \right ], \ \mbox{and} \ q = 1.
\]
The source term $f$ is constructed so that the true solution
is $u = \sin \pi x \sin \pi y.$  We present two experiment sets in which adaptivity 
is driven by a geometric criterion. Namely, the simplices which intersect with the 
quarter circle centered at the origin with radius $0.25$ and $0.05$, in
experiment sets I and II respectively, are repeatedly marked for further refinement.

\noindent $\bullet$ Boundary conditions for the domain in experiment set I:
\begin{eqnarray*}
\Gamma_N & = & \{ (x,y): x = 0, 0 < y < 1 \}  \cup \{ (x,y) : x = 1, 0 < y < 1 \}  \\
\Gamma_D & = & \{ (x,y): 0 \le x \le 1, y = 0 \} \cup \{ (x,y) :  0 \le x \le 1, y  = 1 \}.
\end{eqnarray*}

\noindent  $\bullet$ Boundary conditions for the domain in experiment set II:
\begin{eqnarray*}
\Gamma_N & = & \{ (x,y): 0 \le x \le 1, y = 0 \} \cup \{ (x,y): 0 \le x \le 1, y =1 \} \\
         &   & \cup  \{ (x,y): x = 0, 0 \le y \le 1\} \cup \{ (x,y): x = 1, 0 \le y \le 1\}.
\end{eqnarray*}
\noindent Stopping criterion: $ \| \mbox{error} \|_A < 10^{-7}.$

In experiment set I, red-green refinement subdivides simplices intersecting an arc of 
radius $0.25$ which gives rise to a rapid increase in the number of DOF. Although we have an adaptive refinement strategy, this indeed creates a 
geometric increase  in the number of DOF, see 
Figure~\ref{figure:meshMclabSpecial}. Experiment set II is designed so that a small 
number of DOF is introduced at each level. In order to do 
this, green refinement subdivides simplices intersecting a smaller arc with radius $0.05$.

\begin{table}[htb]
\caption{MCLite iteration counts for various methods, red-green refinement driven by geometric refinement, experiment set I.}
\begin{center}
\begin{tabular}{|l||c|c|c|c|c|c|c|c|}        \hline
\multicolumn{1}{|l||}{Levels}     &
\multicolumn{1}{|c|}{1}    &
\multicolumn{1}{|c|}{2}    &
\multicolumn{1}{|c|}{3}    &
\multicolumn{1}{|c|}{4}    &
\multicolumn{1}{|c|}{5}    &
\multicolumn{1}{|c|}{6}    &
\multicolumn{1}{|c|}{7}    &
\multicolumn{1}{|c|}{8}    \\ \hline \hline
MG       & 1 & 4 & 7 & 7 & 7 & 6 & 6 & 6  \\ \hline
M.BPX    & 1 & 4 & 7 & 7 & 7 & 7 & 6 & 6  \\ \hline
HBMG     & 1 &10 &19 &28 &32 &37 &45 &56  \\ \hline
WMHBMG   & 1 & 6 &12 &13 &16 &17 &17 &17  \\ \hline \hline
PCG-MG     & 1 & 3 & 4 & 5 & 5 & 5 & 5 & 5  \\ \hline
PCG-M.BPX  & 1 & 3 & 5 & 5 & 5 & 5 & 5 & 5  \\ \hline
PCG-HBMG   & 1 & 3 & 7 &10 &12 &14 &15 &16  \\ \hline
PCG-WMHBMG & 1 & 3 & 7 & 7 & 9 & 9 & 9 & 9 \\ \hline \hline
PCG-A.MG & 1 & 8 & 13 & 17 & 20 & 21 & 23 & 24 \\ \hline
PCG-BPX  & 1 & 6 & 12 & 14 & 17 & 17 & 18 & 18 \\ \hline
PCG-HB   & 1 & 5 & 14 & 21 & 26 & 32 & 38 & 41 \\ \hline
PCG-WMHB & 1 & 5 & 12 & 15 & 19 & 20 & 21 & 21 \\ \hline \hline
Nodes & 16 & 19 & 31 & 55 &117 &219 &429 &835  \\ \hline
DOF   &  8 & 10 & 21 & 43 &102 &202 &410 &814  \\ \hline
\end{tabular}
\end{center}
\label{table:ExpSetI}
\end{table}

\begin{table}[ht]
\caption{MCLite iteration counts for various methods, green refinement driven by geometric refinement, experiment set II.}
\begin{center}
\begin{tabular}{|l||c|c|c|c|c|c|c|}        \hline
\multicolumn{1}{|l||}{Levels}     &
\multicolumn{1}{|c|}{1}    &
\multicolumn{1}{|c|}{2}    &
\multicolumn{1}{|c|}{3}    &
\multicolumn{1}{|c|}{4}    &
\multicolumn{1}{|c|}{5}    &
\multicolumn{1}{|c|}{6}    &
\multicolumn{1}{|c|}{7}    \\ \hline
 & 8 & 9 & 10 & 11 & 12 & 13 & 14 \\ \hline \hline
MG    & 1 & 3 & 4 & 3 & 4 & 4 & 3 \\ \hline
      & 4 & 4 & 4 & 4 & 4 & 4 & 4 \\ \hline
M.BPX & 1 & 4 & 4 & 4 & 4 & 4 & 4 \\ \hline
      & 4 & 5 & 5 & 5 & 5 & 5 & 5 \\ \hline
HBMG & 1 & 13 & 14 & 16 & 22 & 25 & 26 \\ \hline
     &30 & 32 & 32 & 36 & 38 & 42 & 44 \\ \hline
WMHBMG  & 1 & 8 &11 &11 &12 &12 &12 \\ \hline
& 13 &15 &15 &15 &15 &15 &15  \\ \hline \hline
PCG-MG    & 1 & 2 & 3 & 3 & 3 & 4 & 3 \\ \hline
        & 3 & 4 & 3 & 4 & 3 & 3 &  3  \\ \hline
PCG-M.BPX & 1 & 2 & 3 & 4 & 4 & 3 & 4 \\ \hline
        & 4 & 4 & 4 & 4 & 4 & 4 &  4  \\ \hline
PCG-HBMG  & 1 & 2 & 5 & 7 & 8 & 9 & 10 \\ \hline
        &10 & 11 & 12 & 11 & 12 & 13 & 13  \\ \hline
PCG-WMHBMG & 1 & 2 & 5 & 6 & 6 & 7 & 7  \\ \hline
         & 8 &  8 & 8 & 8 & 8 & 8 &  8 \\ \hline \hline
PCG-A.MG & 1 & 10 & 13 & 15 & 18 & 20 & 21 \\ \hline
       & 23 & 25 & 26 & 28 & 28 & 28 & 29 \\ \hline
PCG-BPX & 1 & 6 & 10 & 11 & 13 & 14 & 15 \\ \hline 
      & 16 & 18 & 19 & 19 & 20 & 20 & 21 \\ \hline
PCG-HB & 1 & 3 & 9 & 11 & 14 & 18 & 20 \\ \hline
     & 22 & 24 & 27 & 30 & 32 & 34 & 36 \\ \hline
PCG-WMHB & 1 & 3 & 9 & 12 & 14 & 16 & 17  \\ \hline
& 19 & 20 & 20 & 22 & 23 & 23 & 23 \\ \hline \hline
Nodes=DOF & 289 & 290 & 296 & 299 & 309 & 319 & 331 \\ \hline
      & 349 & 388 & 423 & 489 & 567 & 679 & 837 \\ \hline
\end{tabular}
\end{center}
\label{table:ExpSetII}
\end{table}

In all the experiments, we utilize a direct coarsest level solve and smoother 
is a symmetric Gauss-Seidel iteration. The set of DOF on which smoother acts is 
the fundamental difference between the methods. Classical multigrid methods
smooth on all DOF, whereas HB-like methods smooth only on fine DOF. WMHB style
methods smooth as HB methods do. BPX style smoothing is a combination of
multigrid and HB style. 

There are four multiplicative methods under consideration: MG, M.BPX, HBMG, and 
WMHBMG. The following is a guide to the tables and figures below. 
MG will refer to a classical multigrid, in particular corresponds to the
standard V-cycle implementation. HBMG corresponds exactly to the MG algorithm, 
but where pre- and post-smoothing are restricted to fine DOF.
M.BPX refers to multiplicative version of BPX. Smoother is 
restricted to fine DOF and their immediate coarse neighbors which are often 
called as the {\em 1-ring} neighbors. 1-ring neighbors of the fine nodes can be 
directly determined by the sparsity pattern of the fine-fine subblock $A_{22}$
of the stiffness matrix. The set of DOF over which BPX method smooths is simply 
the union of the column locations of nonzero entries corresponding to fine DOF.
Using this observation, HBMG smoother can easily be modified to be a BPX 
smoother. WMHBMG is similar to HBMG, in that both are multiplicative methods
in the sense of Definition~\ref{defn:Bj}, but the difference is in the basis 
used. In particular, the change of basis matrices are different as a result
of the wavelet stabilization, where $L_2$-projection to coarser finite
element spaces is approximated by two Jacobi iterations.

\begin{figure}[h]
\begin{center}
\mbox{
\myfigpdf{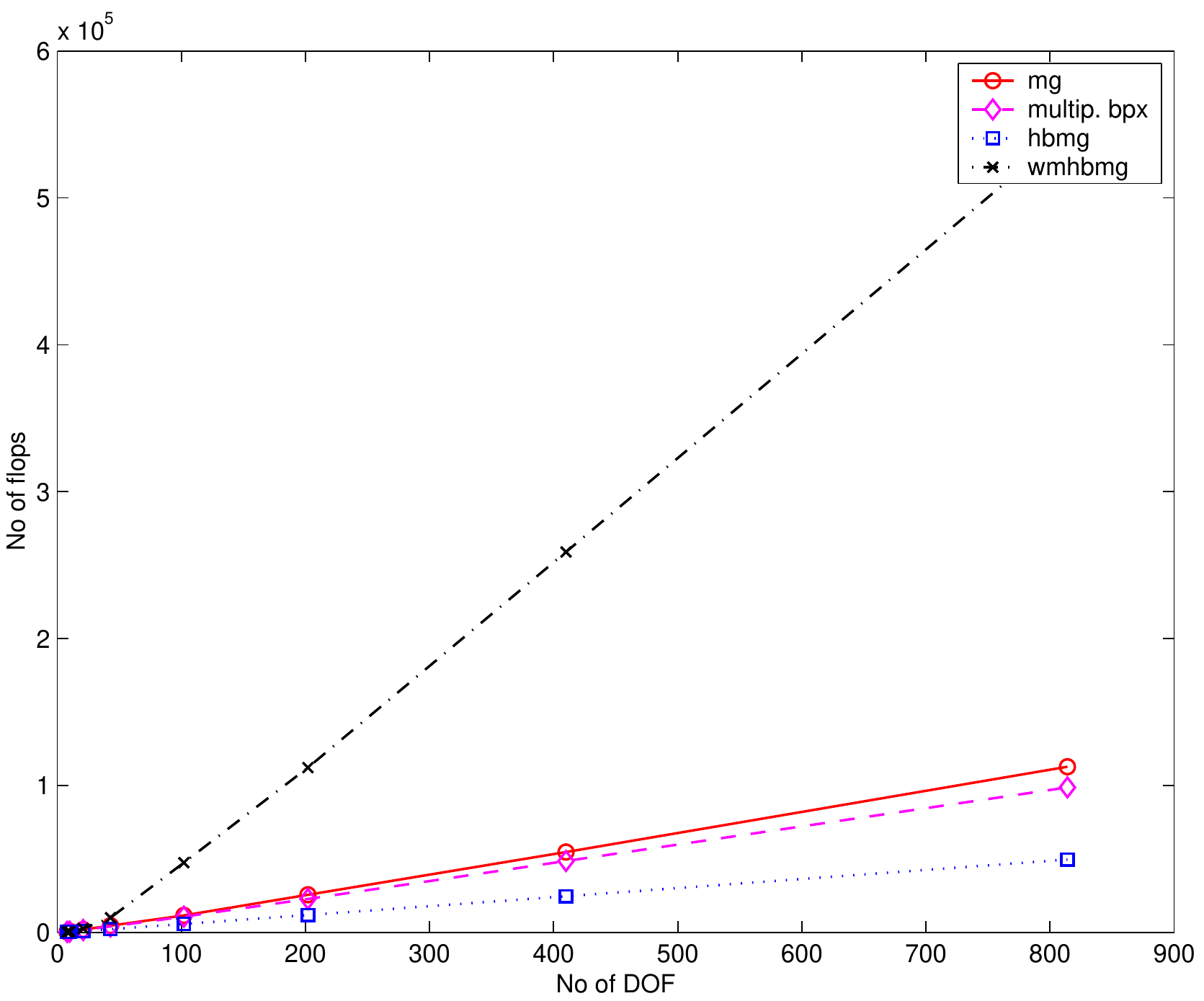}{6cm} 
\myfigpdf{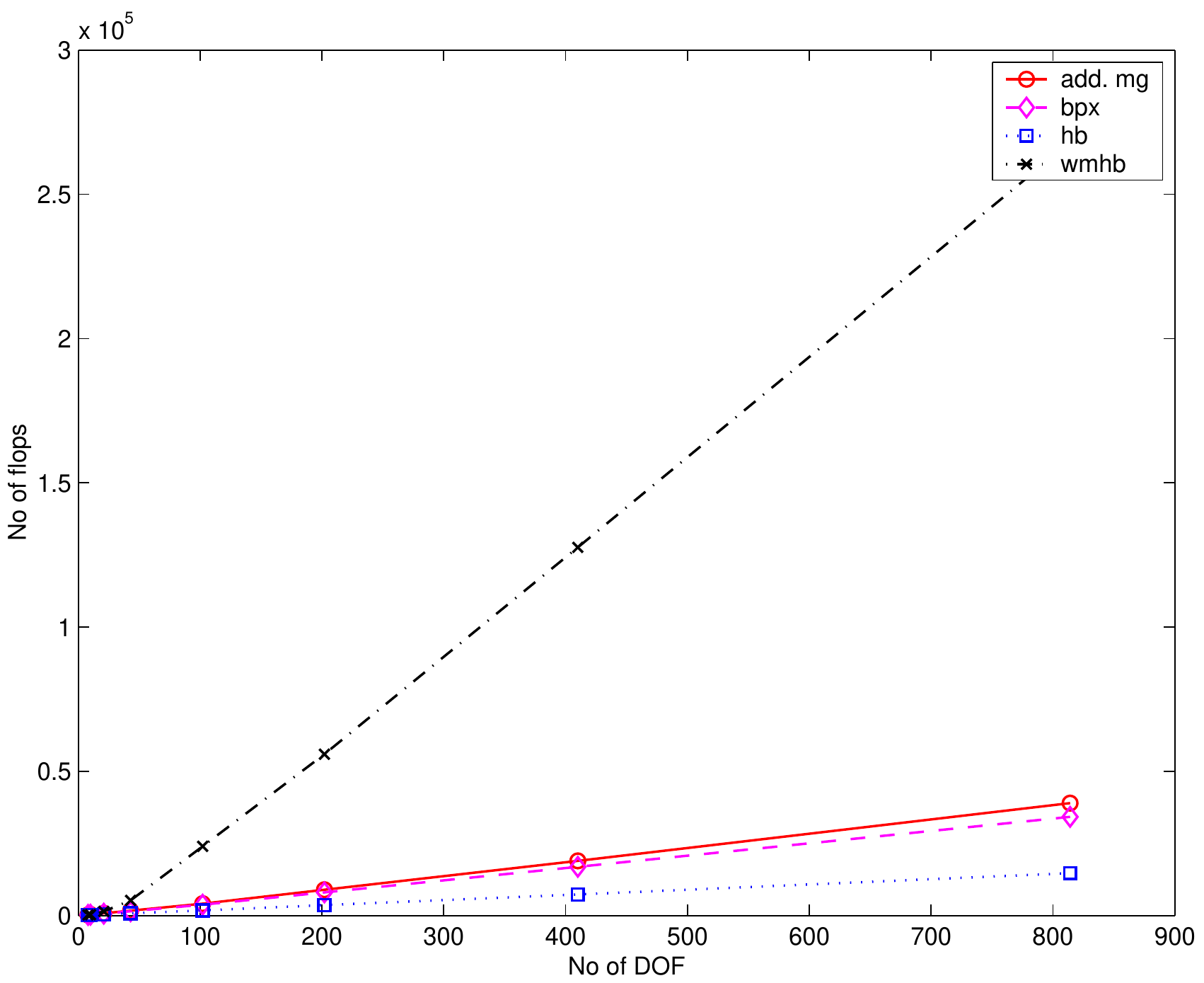}{6cm} 
}
\caption{Flop counts for single iteration of multiplicative (left) and additive (right) methods, experiment set I.}
\label{figure:flopCountMulAddExpI}
\end{center}
\end{figure}

\begin{figure}[ht]
\begin{center}
\mbox{
\myfigpdf{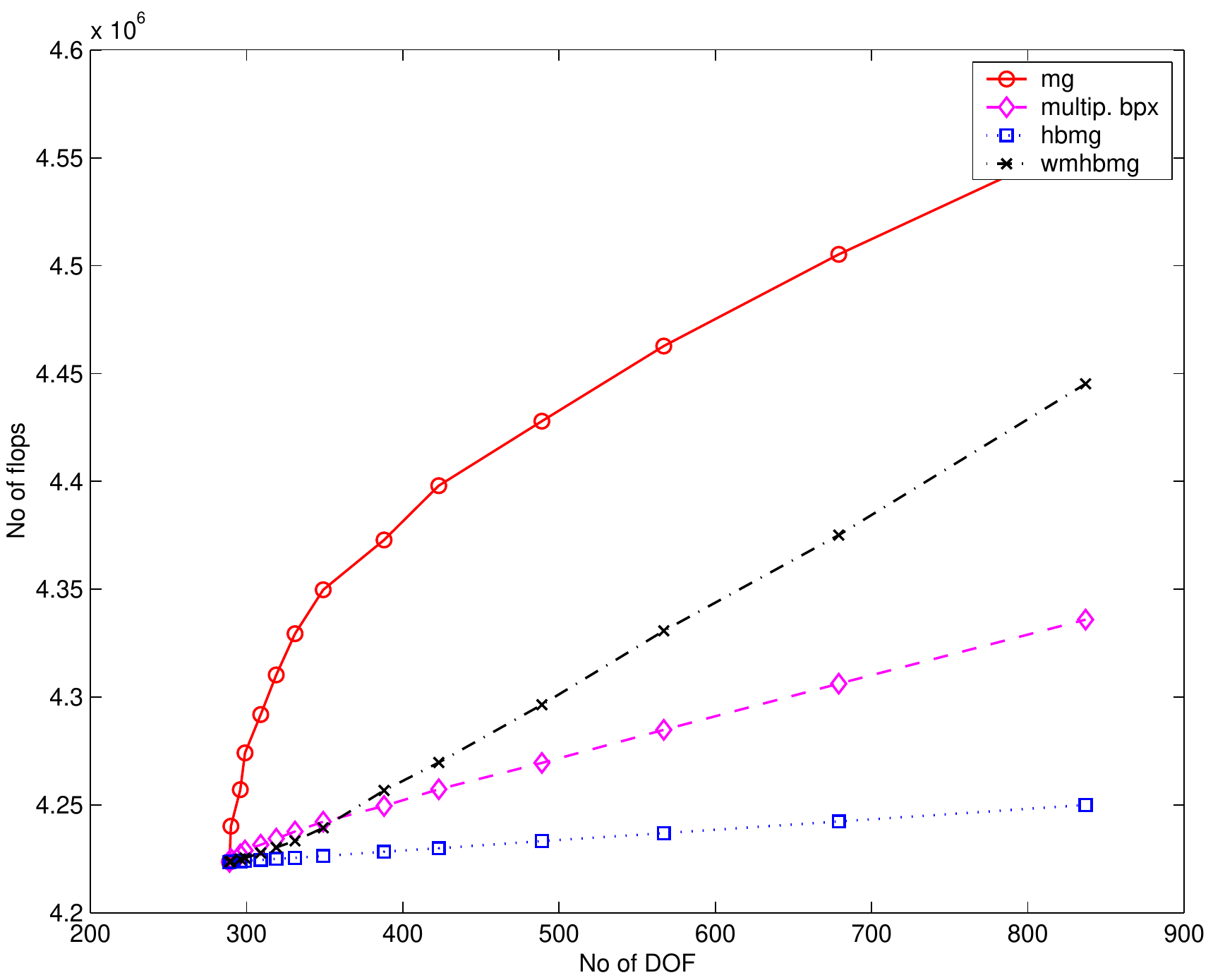}{6cm} 
\myfigpdf{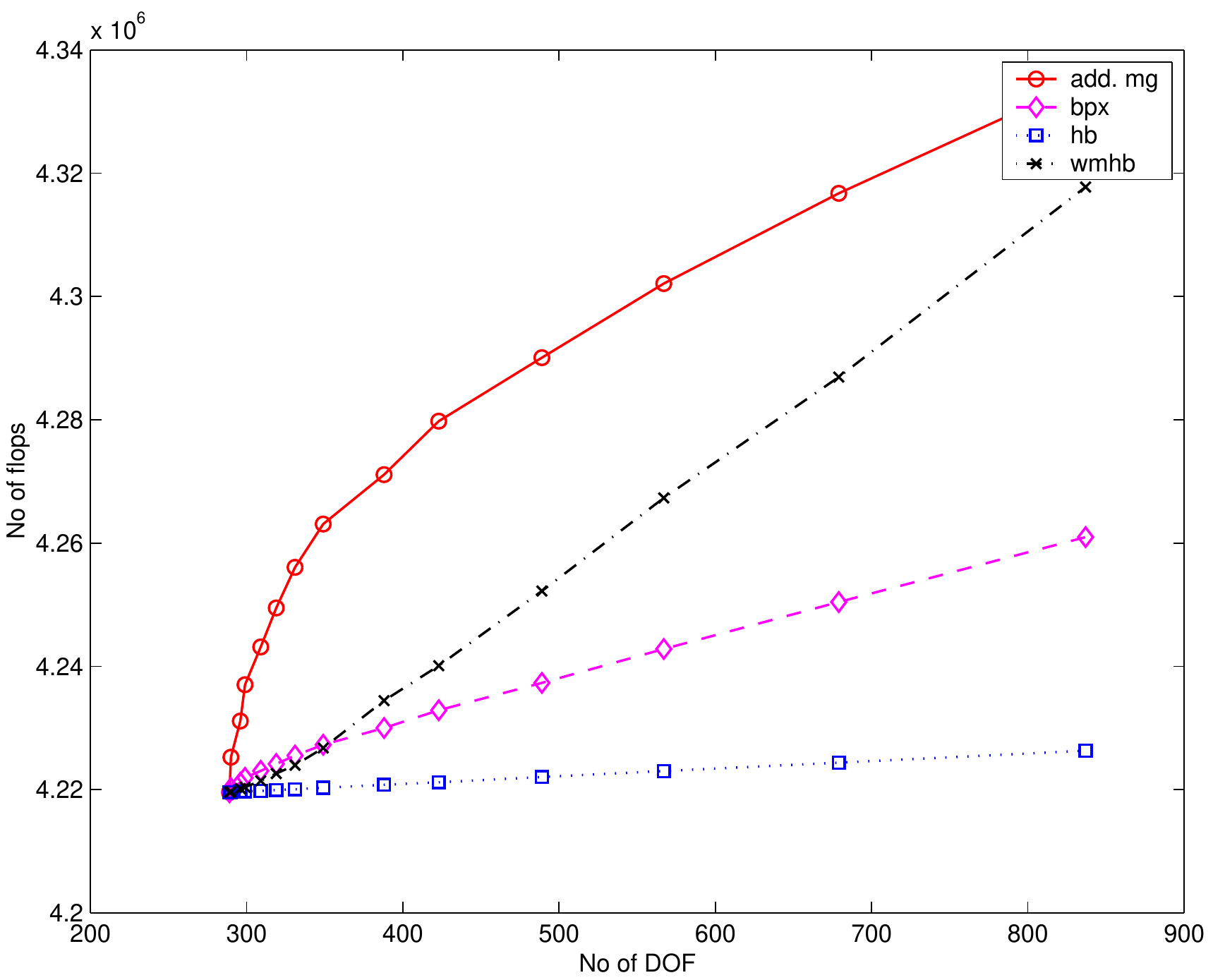}{6cm} 
}
\caption{Flop counts for single iteration of multiplicative (left) and additive (right) methods, experiment set II.}
\label{figure:flopCountMulAddExpII}
\end{center}
\end{figure}

PCG stands for the preconditioned conjugate gradient method.
PCG-A.MG, PCG-BPX, PCG-HB, and PCG-WMHB involve the use of additive MG, PBX, 
HB, and WMHB as preconditioners for CG, respectively. In the sense of 
Definition~\ref{defn:Dj}, HB and WMHB are additive versions of HBMG and WMHBMG 
respectively. Each preconditioner is implemented in a manor similar to that
described in~\cite{VaWa2,VaWa3}.

Finally, note that {\em Nodes} denotes the total number of nodes in the
simplicial mesh, including Dirichlet and Neumann nodes.
The iterative methods view DOF as the union of the unknowns
corresponding to interior and Neumann/Robin boundary DOF,
and these are denoted as such.

The refinement procedure utilized in the experiments are fundamentally the same as 
the 2D red-green described in~\cite{AkHo02-p1,AkHo02-p2}. We, however, remove
the restrictive conditions that the simplices for level $j+1$ have to be created
from the simplices at level $j$ and the bisected (green refined) simplices cannot 
be further refined. Even in this case the claimed results seem to hold.
Experiments are done in MCLite module of the FEtk package.
Several key routines from the MCLite software, used to produce most of
the numerical results in this paper, are given in the appendix.

Iteration counts are reported in Tables~\ref{table:ExpSetI}
and~\ref{table:ExpSetII}. The optimality of M.BPX, BPX, WMHBMG and WMHB is 
evidenced in each of the experiments with the constant number of iterations, 
independent of the number of DOF. HB and HBMG methods suffer from
a logarithmic increase in the number of iterations.
Among all the methods tested, the M.BPX is the closest to MG in terms of
low iteration counts.

\begin{figure}[ht]
\begin{center}
\mbox{
\myfigpdf{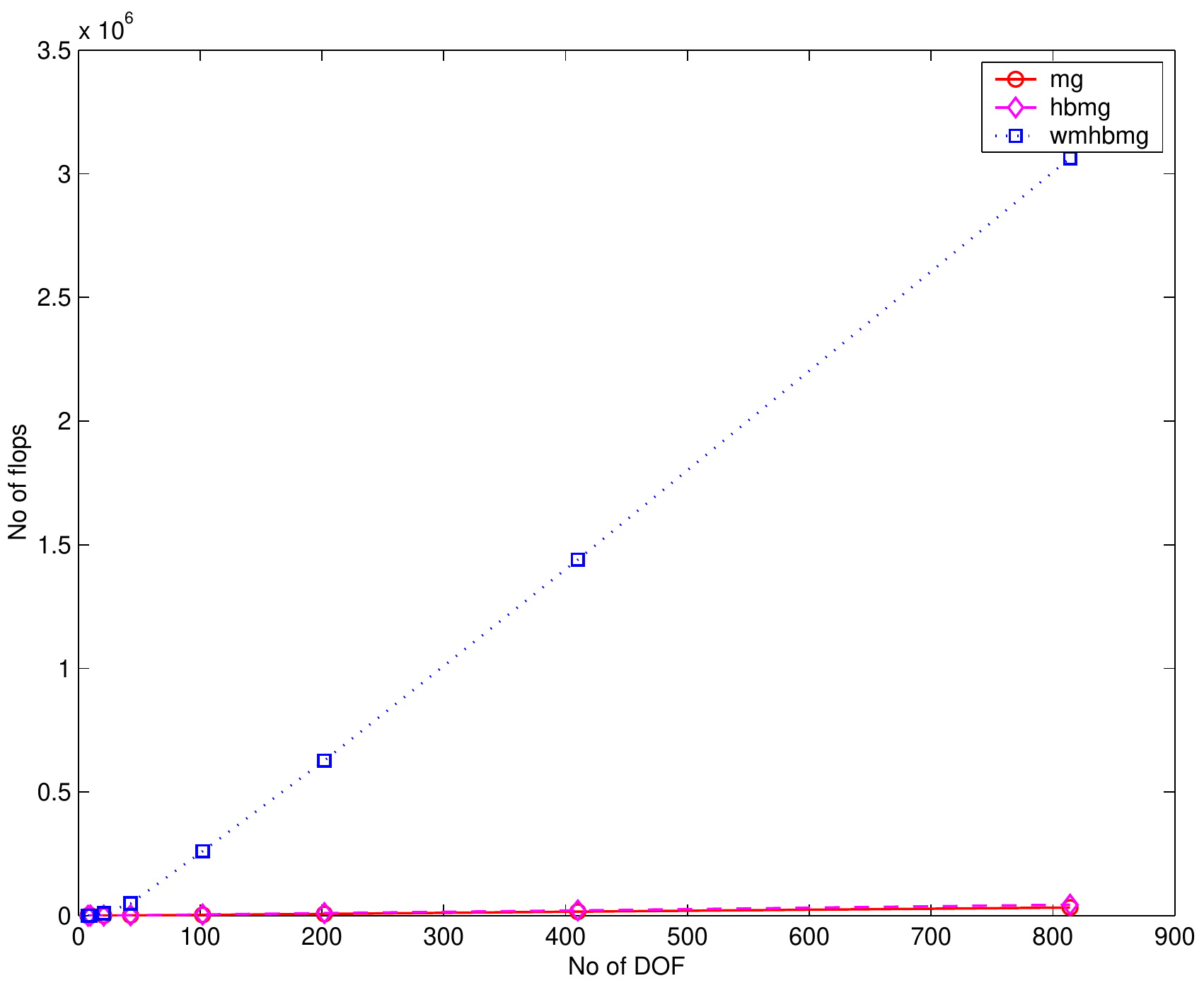}{6cm} 
\myfigpdf{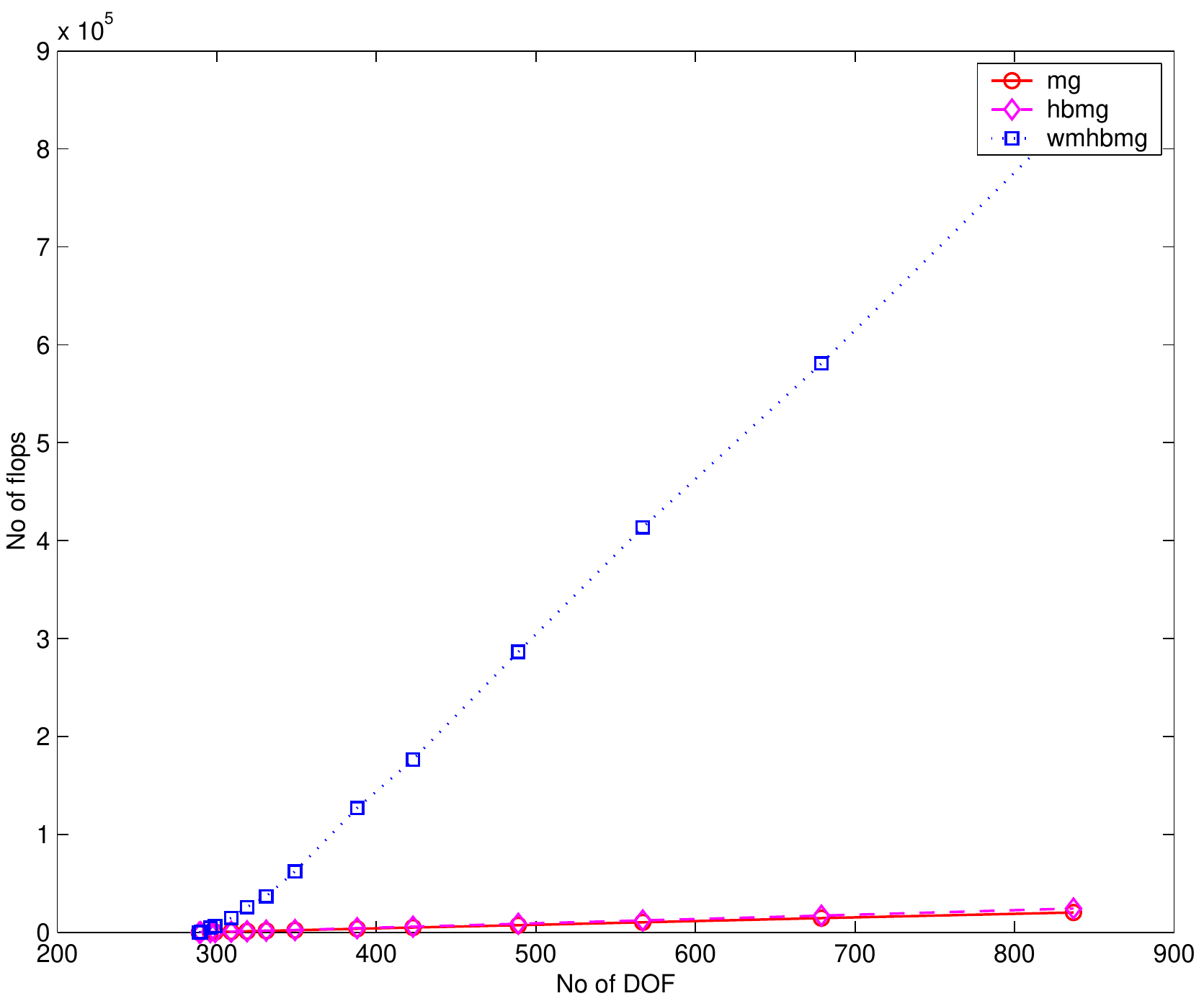}{6cm} 
}
\caption{Flop counts for variational conditions for experiment set I (left) and experiment set II (right).}
\label{figure:variationCountIandII}
\end{center}
\end{figure}

\begin{figure}[ht]
\begin{center}
\mbox{
\myfigpdf{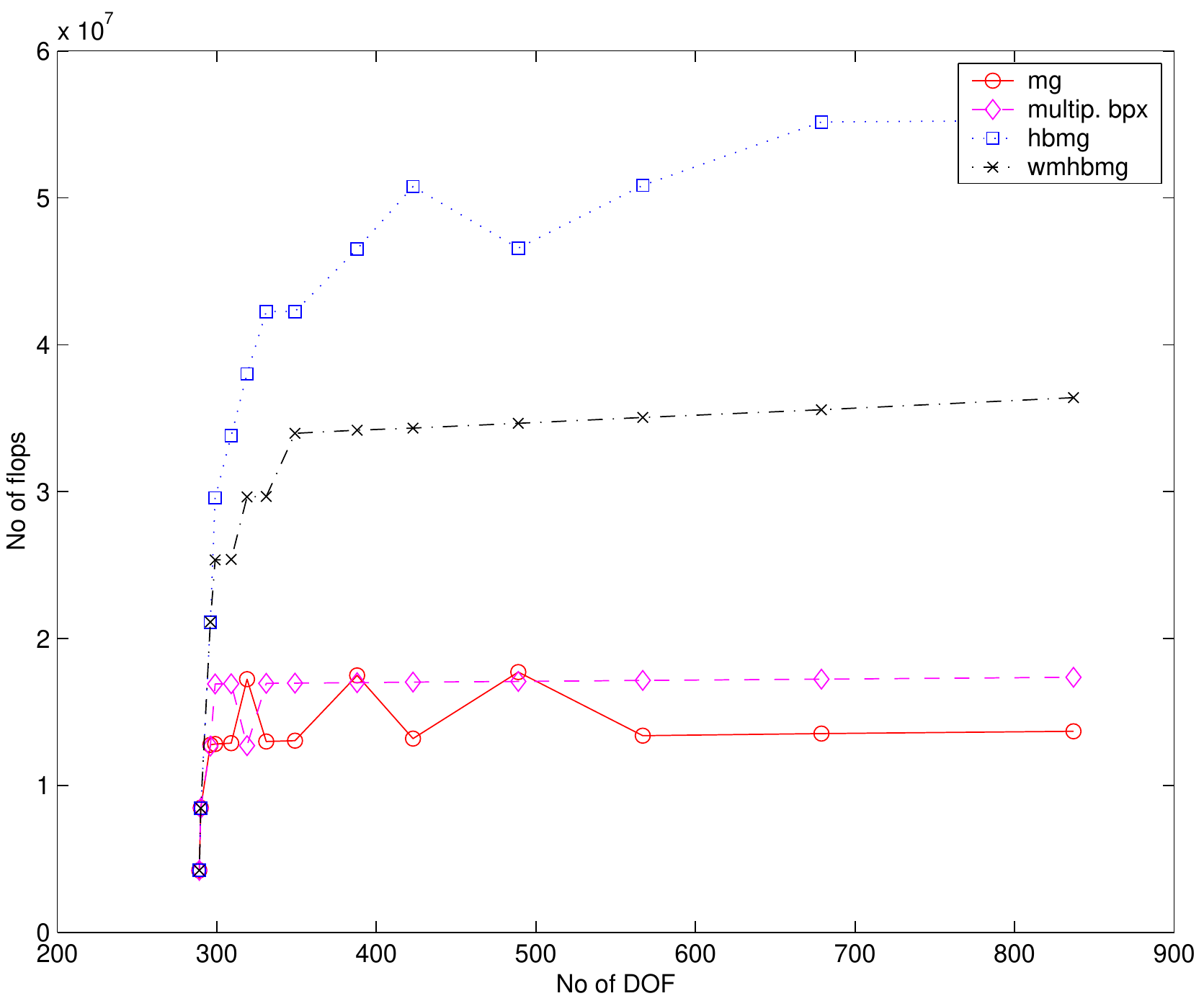}{6cm} 
\myfigpdf{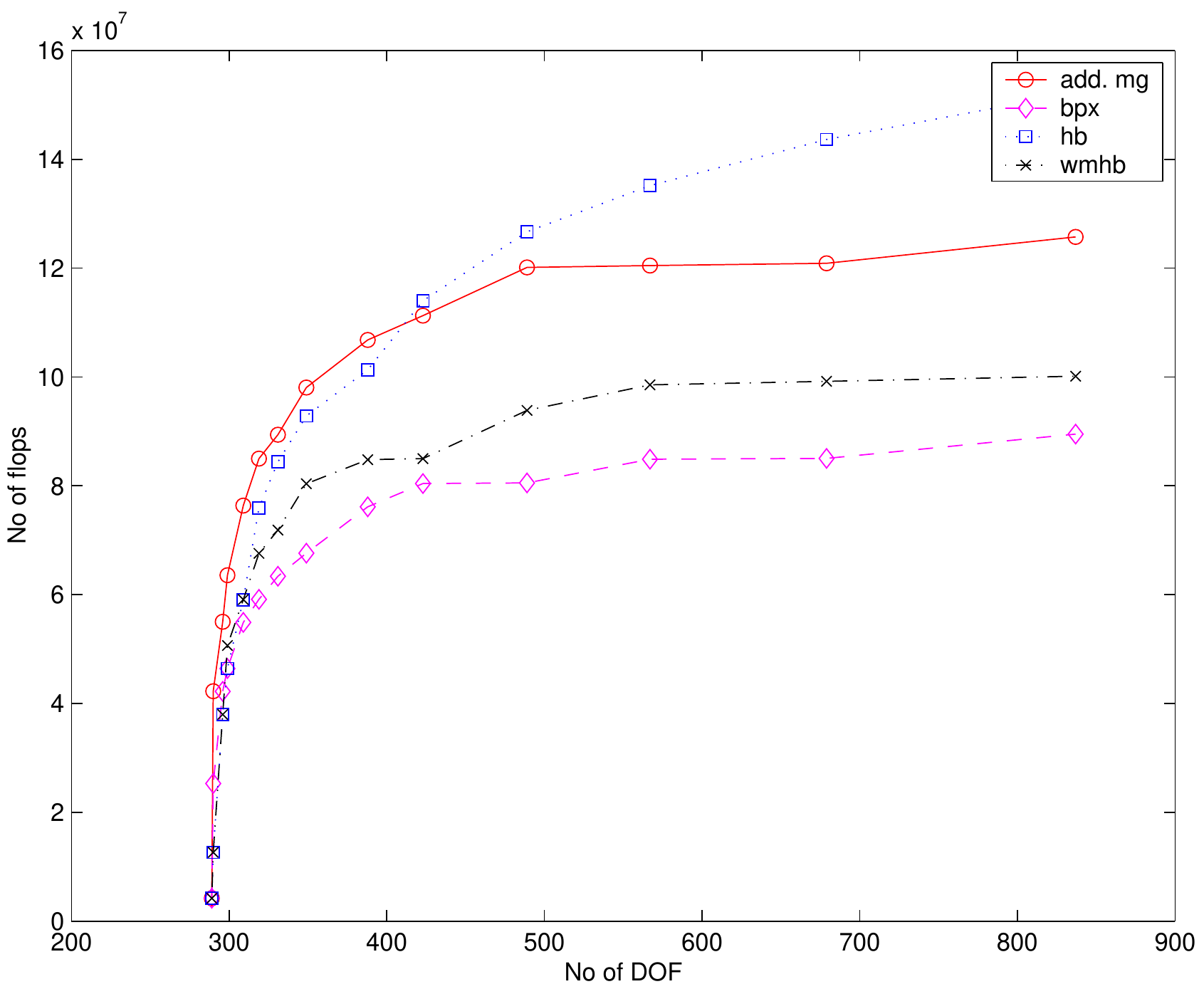}{6cm} 
}
\caption{Total flop counts of  preconditioned multiplicative (left) and additive (right) methods, experiment set II.}
\label{figure:TotalFlopCountMulAddExpII}
\end{center}
\end{figure}

However, it should be clearly noted that in the experiments we present below, 
the cost per iteration of the various methods can differ substantially. We
report flop counts of a single iteration of the above methods, see 
Figures~\ref{figure:flopCountMulAddExpI} and~\ref{figure:flopCountMulAddExpII}.
In experiment set I, the cost per iteration is linear for all the methods.
The WMHB and WMHBMG methods are the most expensive ones. 
We would like to emphasize that the refinement in experiment set I cannot be a 
good example for adaptive refinement given the geometric increase in the number
of DOF. MG exploits this geometric increase and enjoys a linear computational 
complexity. Experiment set II is more realistic in the sense that the 
refinement is highly adaptive and introduces a small number of DOF at each 
level. One can now observe a suboptimal (logarithmic) computational complexity 
for MG-like methods in such realistic scenarios. In accordance with the 
theoretical justification, under highly adaptive refinement MG methods will 
asymptotically be suboptimal. Moreover, storage complexity prohibitively 
prevents MG-like methods from being a viable tool for large and highly adaptive
settings.

Coarser representations of the finest level system~(\ref{ALGproblem})
are algebraically formed by enforcing variational conditions. Some methods 
require further stabilizations in the from of matrix-matrix products. These form
the so-called {\em preprocessing} step in multilevel methods. The computational 
cost of variational conditions is the same regardless of having a 
multiplicative or an additive version of the same method. This computational
cost is orders of magnitude cheaper than the cost of a single iteration. 
However, this is the 
step where the storage complexity can dominate the overall complexity. Due to 
memory bandwidth problems on conventional machines, one should be very careful
with the choice of datastructures.  Since only the $A_{11}=A_{coarse}$ subblock
of $A$ is formed for the next coarser level, the cost of variational 
conditional for MG, M.BPX, A.MG, and BPX is the cheapest among all the methods.
On the other hand, HBMG and HB require stabilizations of $A_{12}$ and 
$A_{21}$ using to the hierarchical basis. The WMHBMG and WMHB methods are more 
demanding by requiring stabilizations of $A_{12}$, $A_{21}$, and $A_{22}$ 
using the wavelet modified hierarchical basis. Wavelet structure creates denser
change of basis matrix than that of the hierarchical basis. Therefore, 
preprocessing in the WMHB and WMHBMG methods is the most expensive among all 
the methods.

\section{Conclusion}
  \label{sec:conc}

In this paper, we examined a number of additive and multiplicative
multilevel iterative methods and preconditioners in the setting of
two-dimensional local mesh refinement. While standard multilevel methods are 
effective for uniform refinement-based discretizations of elliptic equations,
they tend to be less effective for algebraic systems which arise from 
discretizations on locally refined meshes, losing both their optimal behavior
in both storage or computational complexity. Our primary focus here was on 
BPX-style additive and multiplicative multilevel preconditioners,
and on various stabilizations of the additive and multiplicative
hierarchical basis method, and their use in the local mesh refinement setting.
In the first two papers of this trilogy, it was shown that both BPX and
wavelet stabilizations of HB have uniformly bounded conditions numbers on
several
classes of locally refined 2D and 3D meshes based on fairly standard
(and easily implementable) red and red-green mesh refinement algorithms.
In this third article of the trilogy, we described in detail the
implementation of these types
of algorithms, including detailed discussions of the datastructures and
traversal algorithms we employ for obtaining optimal storage and computational
complexity in our implementations.
We showed how each of the algorithms can be implemented using standard
datatypes available in languages such as C and FORTRAN, so that the
resulting algorithms
have optimal (linear) storage requirements, and so that the resulting
multilevel method or preconditioner can be applied with optimal (linear)
computational costs.

Our implementations were performed in both C and MATLAB using the
Finite Element ToolKit (FEtk), an open source finite element software
package. We presented a sequence of numerical experiments illustrating
the effectiveness of a number of BPX and stabilized HB variants for several
examples requiring local refinement. As expected, multigrid methods most 
effective in terms of iteration counts remaining constant as the DOF
increase, but the suboptimal complexity per iteration in the local refinement 
setting makes the BPX methods the most attractive.
In addition, both the additive and multiplicative WMHB-based methods
and preconditioners demonstrated similar constant iteration requirements
with increasing DOF, yet the cost per iteration remains
optimal (linear) even in the local refinement setting.
Consequently in highly adaptive regimes, the BPX methods prove to be the 
most effective, and the WMHB methods become the second best effective.

\begin{center}
\begin{figure}[ht]
\mbox{
\begin{minipage}{15.0cm}
\begin{center}
\mbox{\myfigpdf{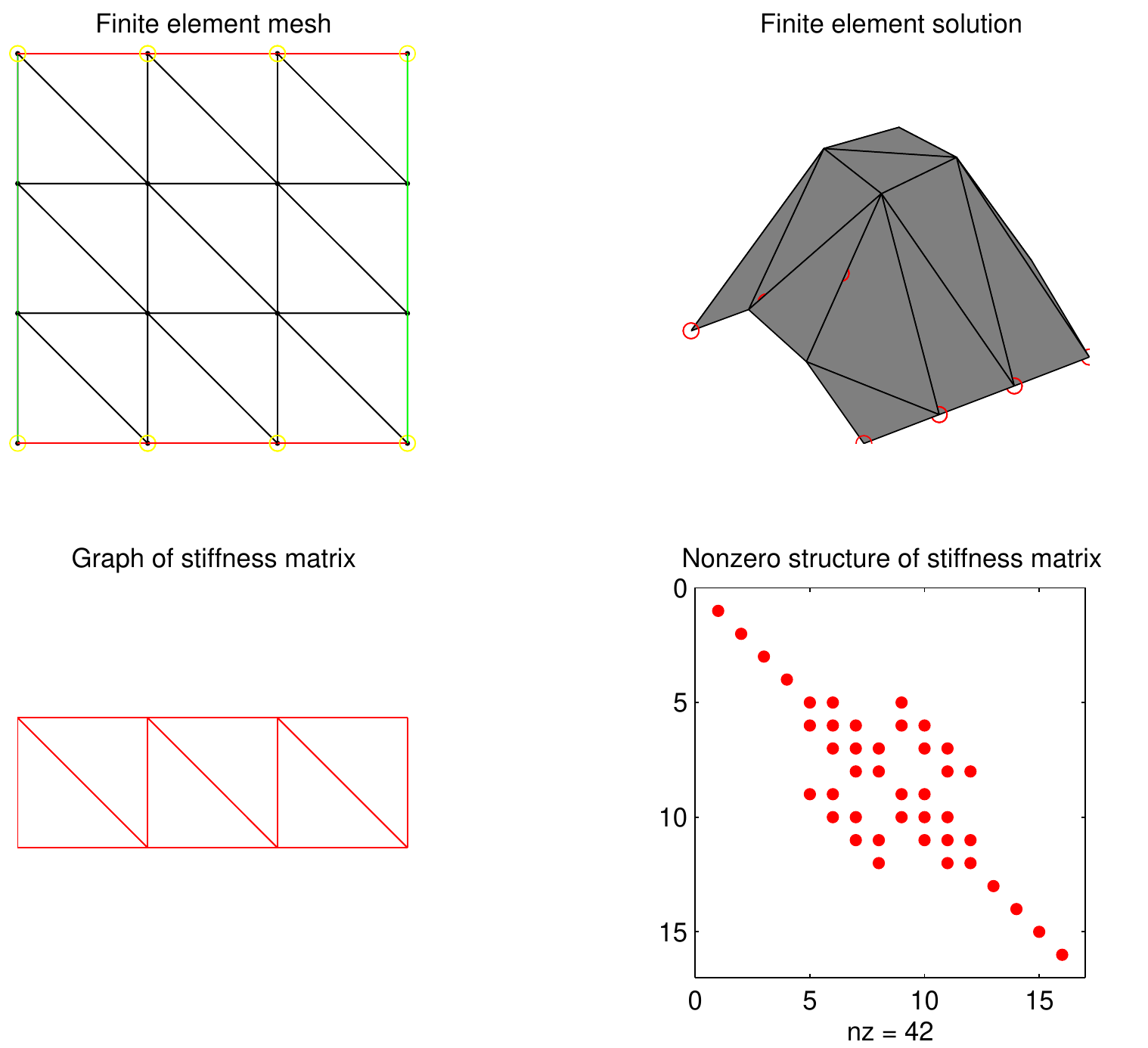}{6cm}} \hspace{1cm} 
\mbox{\myfigpdf{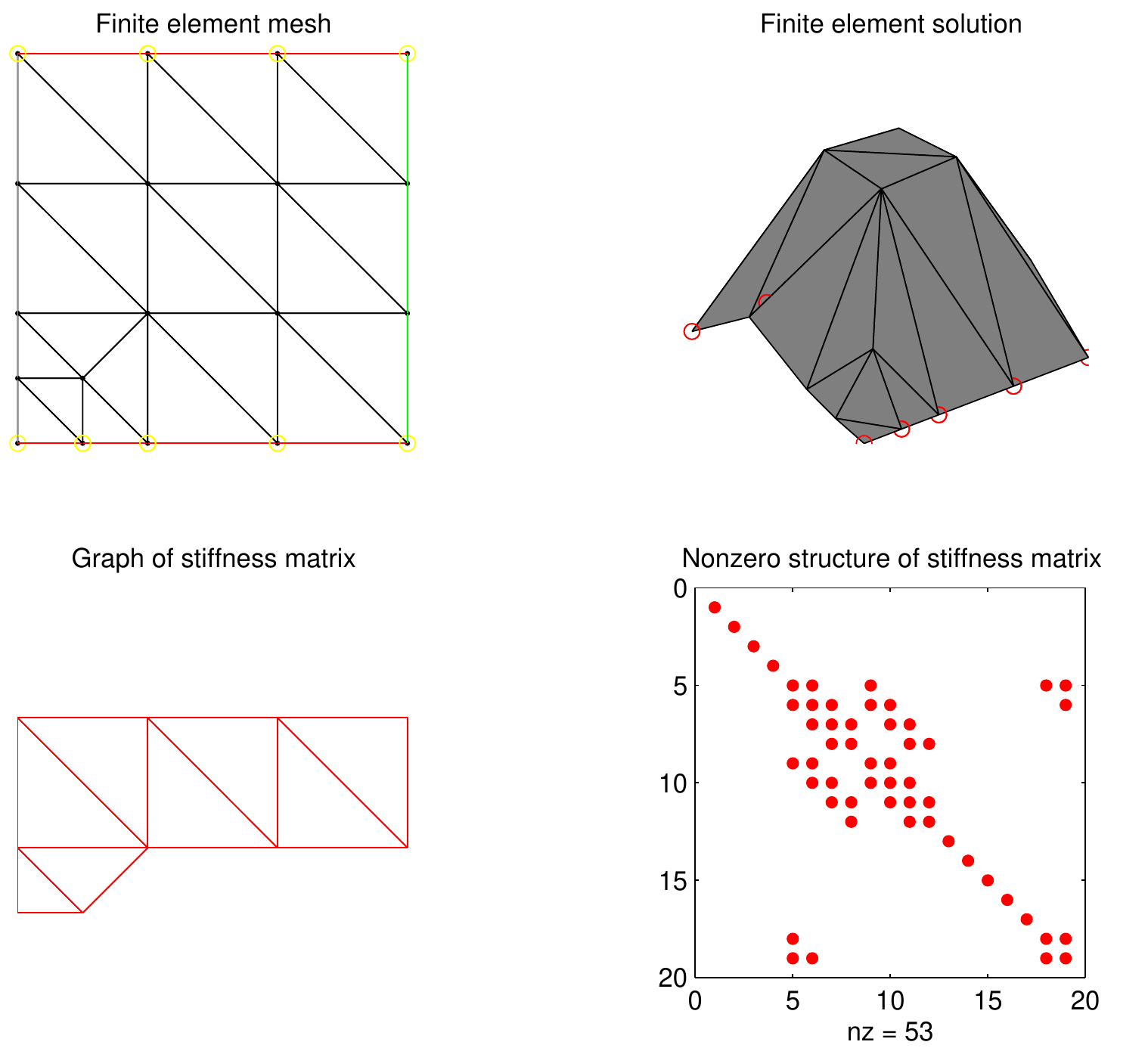}{6cm}} \\ 
\vspace{.1cm}
\mbox{\myfigpdf{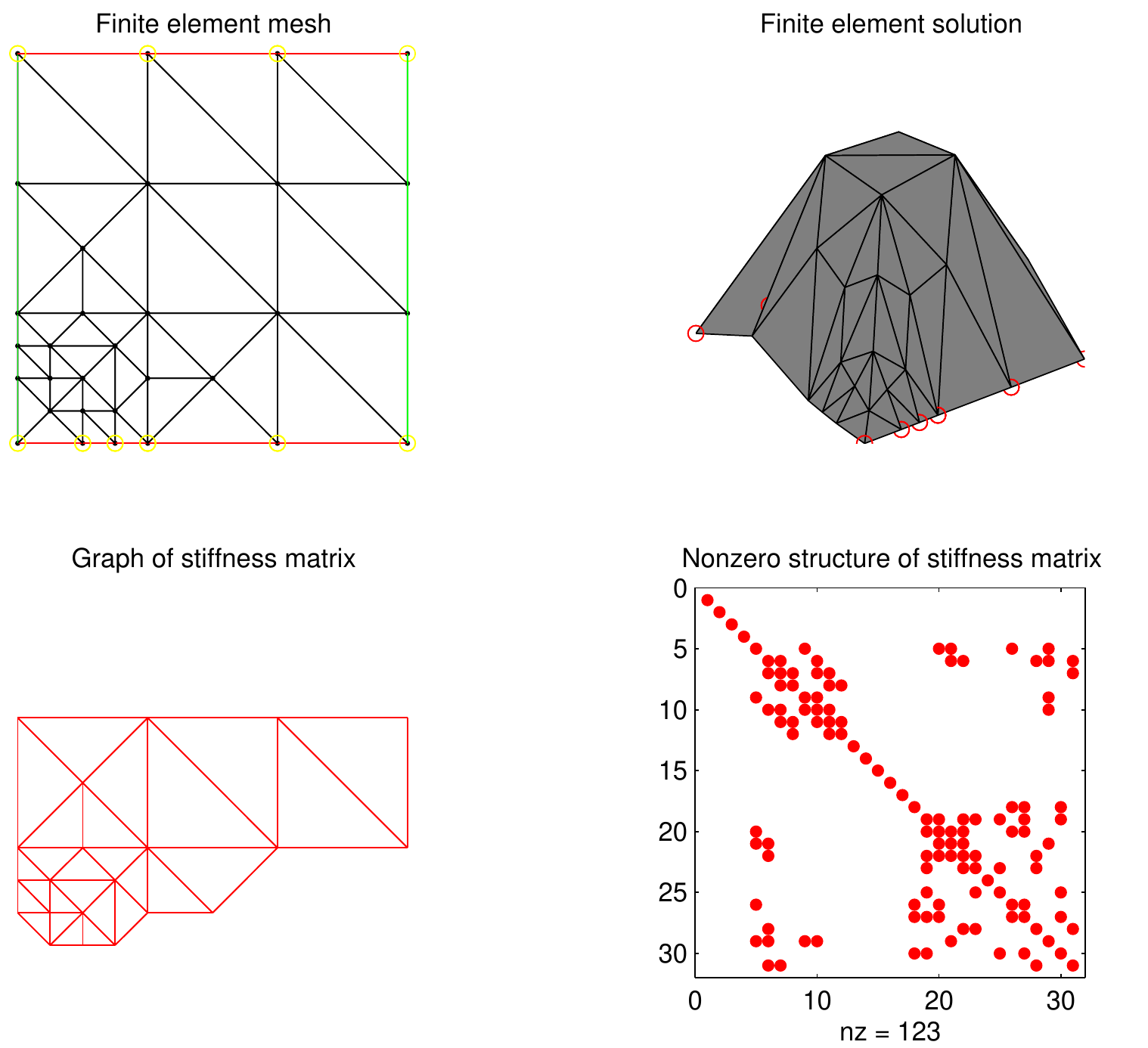}{12cm}} 
\end{center}
\end{minipage}
}
\caption{Dialog box from MCLite for experiment set I.} 
\label{figure:charts}
\end{figure}
\end{center}

\section*{Acknowledgments}
The authors thank R.~Bank for many enlightening discussions.

\bibliographystyle{siam}
\bibliography{m}

\begin{thebibliography}{10}

\bibitem{Ak01phd}
{\sc B.~Aksoylu}, {\em Adaptive Multilevel Numerical Methods with Applications
  in Diffusive Biomolecular Reactions}, PhD thesis, Department of Mathematics,
  University of California, San Diego, La Jolla, CA, 2001.

\bibitem{AkHo02-p1}
{\sc B.~Aksoylu and M.~Holst}, {\em {An odyssey into local refinement and
  multilevel preconditioning I: Optimality of the BPX preconditioner}}, SIAM J.
  Numer. Anal.,  (2002).
\newblock in review.

\bibitem{AkHo02-p2}
\leavevmode\vrule height 2pt depth -1.6pt width 23pt, {\em {An odyssey into
  local refinement and multilevel preconditioning II: Stabilizing hierarchical
  basis methods}}, SIAM J. Numer. Anal.,  (2002).
\newblock in review.

\bibitem{BaDu80}
{\sc R.~E. Bank and T.~Dupont}, {\em Analysis of a two-level scheme for solving
  finite element equations}, tech. report, Center for Numerical Analysis,
  University of Texas at Austin, 1980.
\newblock CNA--159.

\bibitem{BaDuYs88}
{\sc R.~E. Bank, T.~Dupont, and H.~Yserentant}, {\em The hierarchical basis
  multigrid method}, Numer. Math., 52 (1988), pp.~427--458.

\bibitem{BoYs93}
{\sc F.~Bornemann and H.~Yserentant}, {\em A basic norm equivalence for the
  theory of multilevel methods}, Numer. Math., 64 (1993), pp.~455--476.

\bibitem{BrPa93}
{\sc J.~H. Bramble and J.~E. Pasciak}, {\em New estimates for multilevel
  algorithms including the {V}-cycle}, Math. Comp., 60 (1993), pp.~447--471.

\bibitem{BPX90}
{\sc J.~H. Bramble, J.~E. Pasciak, and J.~Xu}, {\em Parallel multilevel
  preconditioners}, Math. Comp., 55 (1990), pp.~1--22.

\bibitem{DaKu92}
{\sc W.~Dahmen and A.~Kunoth}, {\em Multilevel preconditioning}, Numer. Math.,
  63 (1992), pp.~315--344.

\bibitem{DeEiGiLiLiu99}
{\sc J.~W. Demmel, S.~C. Eisenstat, J.~R. Gilbert, X.~S. Li, and J.~W.~H. Liu},
  {\em A supernodal approach to sparse partial pivoting}, SIAM J. Matrix Anal.
  Appl., 20 (1999), pp.~720--755.

\bibitem{DuGrLe89}
{\sc I.~S. Duff, R.~G. Grimes, and J.~G. Lewis}, {\em Sparse matrix test
  problems}, ACM Trans. Math. Softw., 15 (1989), pp.~1--14.

\bibitem{Hols2001a}
{\sc M.~Holst}, {\em Adaptive numerical treatment of elliptic systems on
  manifolds}, Advances in Computational Mathematics, 15 (2001), pp.~139--191.

\bibitem{Ja92}
{\sc S.~Jaffard}, {\em Wavelet methods for fast resolution of elliptic
  problems}, SIAM J. Numer. Anal., 29 (1992), pp.~965--986.

\bibitem{Os94book}
{\sc P.~Oswald}, {\em Multilevel Finite Element Approximation Theory and
  Applications}, Teubner Skripten zur Numerik, B. G. Teubner, Stuttgart, 1994.

\bibitem{St95}
{\sc R.~Stevenson}, {\em Robustness of the additive multiplicative frequency
  decomposition multi-level method}, Computing, 54 (1995), pp.~331--346.

\bibitem{St97}
\leavevmode\vrule height 2pt depth -1.6pt width 23pt, {\em A robust
  hierarchical basis preconditioner on general meshes}, Numer. Math., 78
  (1997), pp.~269--303.

\bibitem{VaWa2}
{\sc P.~S. Vassilevski and J.~Wang}, {\em Stabilizing the hierarchical basis by
  approximate wavelets, {I}: Theory}, Numer. Linear Alg. Appl., 4 Number 2
  (1997), pp.~103--126.

\bibitem{VaWa1}
\leavevmode\vrule height 2pt depth -1.6pt width 23pt, {\em Wavelet-like methods
  in the design of efficient multilevel preconditioners for elliptic {PDEs}},
  in Multiscale Wavelet Methods For Partial Differential Equations, W.~Dahmen,
  A.~Kurdila, and P.~Oswald, eds., Academic Press, 1997, ch.~1, pp.~59--105.

\bibitem{VaWa3}
\leavevmode\vrule height 2pt depth -1.6pt width 23pt, {\em Stabilizing the
  hierarchical basis by approximate wavelets, {II}: Implementation and
  numerical experiments}, SIAM J. Sci. Comput., 20 Number 2 (1998),
  pp.~490--514.

\bibitem{XuQi94}
{\sc J.~Xu and J.~Qin}, {\em Some remarks on a multigrid preconditioner}, SIAM
  J. Sci. Comput., 15 (1994), pp.~172--184.

\bibitem{Ys86}
{\sc H.~Yserentant}, {\em On the multilevel splitting of finite element
  spaces}, Numer. Math., 49 (1986), pp.~379--412.

\end{thebibliography}


\vspace*{0.5cm}

\end{document}